\documentclass[a4paper]{article}
\usepackage[affil-it]{authblk}

\usepackage[english]{babel} 
\usepackage{amsmath}
\usepackage{amsfonts}
\usepackage{amssymb}
\usepackage{mathrsfs}
\usepackage{graphicx}
\usepackage{enumerate}
\usepackage{xcolor}
\usepackage{hyperref}
\usepackage[numbers,square]{natbib}
\usepackage{amsthm}
\usepackage{faktor}
\usepackage{dirtytalk}
\usepackage{fancyhdr}

\usepackage[Sonny]{fncychap}

\usepackage{tikz}

\usepackage{mathtools,extarrows}
\usepackage{titlesec} 

\titleformat{\chapter}[display]
  {\bfseries\Large}
  {\filright\MakeUppercase{\chaptertitlename} \Huge\thechapter}
{1ex}
  {\titlerule\vspace{1ex}\filleft}  
  [\vspace{1ex}\titlerule]
\DeclarePairedDelimiter\ceil{\lceil}{\rceil}
\DeclarePairedDelimiter\floor{\lfloor}{\rfloor}

\newcommand*\circled[1]{\tikz[baseline=(char.base)]{
            \node[shape=circle,draw,inner sep=1pt] (char) {#1};}}

\newenvironment{pf}{Proof:}{$\hspace*{\fill}\Diamond$}
\newtheorem{df}{Definition}[section]
\newtheorem{prop}{Proposition}[section]
\newtheorem{cor}{Corollary}[section]
\newtheorem{ex}{Example}[section]
\newtheorem{rmk}{Remark}[section]
\newtheorem{algo}{Algorithm}[section]
\newtheorem{thm}{Theorem}[section]
\newtheorem{lem}{Lemma}[section]

\begin{document}

\title{Explicit description of isogeny and isomorphism classes of Drinfeld modules over finite field}

%\author{Yossi Farjoun%
%  \thanks{Electronic address: \texttt{yfarjoun@math.mit.edu}; Corresponding author}}
%\affil{G. Mill\'an Institute of Fluid Dynamics,\\ Nanoscience and Industrial
%Mathematics,\\ Universidad Carlos III de Madrid, Spain}

\author{Sedric Nkotto Nkung Assong
  \thanks{Electronic address: \texttt{sedric.assong@mathematik.uni-kassel.de}}}
\affil{Institute of Mathematics, University of Kassel, Germany}

\date{Dated: \today}

\maketitle

\begin{abstract}

When travelling from the number fields theory to the function fields theory, one cannot miss the deep analogy between rank 1 Drinfeld modules and the group of root of unity and the analogy between rank 2 Drinfeld modules and elliptic curves. But so far, there is no known structure in number fields theory that is analogous to the Drinfeld modules of higher rank $r \geq 3$.\\
In this paper we investigate the classes of those Drinfeld modules of higher rank $r \geq 3$. We describe explicitly the Weil polynomials defining the isogeny classes of rank $r$ Drinfeld modules for any rank $r \geq 3$. our explicit description of the Weil polynomials depends heavily on Yu's classification of isogeny classes (analogue of Honda-Tate at abelian varieties). Actually Yu has also explicitly did that work for $r=2$. \\
To complete the classification, we define the new notion of fine isomorphy invariants for any rank $r$ Drinfeld module and we prove that the fine isomorphy invariants together with J-invariants completely determine the $L$-isomorphism classes of rank $r$ Drinfeld modules defined over the finite field $L$.

\end{abstract}

\section{Isogeny classes}
\label{first section}
\textbf{Notations}:\\

$
\begin{array}{rl}
A=\mathbb{F}_q[T]: & \text{ Ring of univariate polynomials in $T$ over a finite field } \mathbb{F}_q=\mathbb{F}_{p^*},~p~ prime.\\
k=\mathbb{F}_q(T): & \text{ Rational function field over } \mathbb{F}_q.\\
L: & \text{ Finite $A$-field}.\\
\mathfrak{p}_v: & \text{ Kernel of the $\mathbb{F}_q$-algebra homomorphism $\gamma$ defining the $A$-field $L$}.\\
~ & \text{The same notation is used for the ideal and its generator.}\\
v: & \text{the place of $k$ defined by $\mathfrak{p}_v$}\\
\infty: & \text{The place at infinity of $k$}.\\
m : & \text{The degree of $L$ over $A/\mathfrak{p}_v$ i.e. } m=[L : A/\mathfrak{p}_v]   

\end{array}
$

\begin{df}~\\
Let $L\{\tau\}$ be the ring of $\mathbb{F}_q$-linear polynomials spanned by $\{ \tau^i,~ i=1, 2, \cdots \}$ where $\tau$ is the polynomial defined by $\tau(x) = x^q$.
\begin{itemize}
\item A Drinfeld module $\phi$ over $L$ is an $\mathbb{F}_q$-algebra homomorphism\\ $\phi : A \longrightarrow L\{\tau \}$ such that\\
- $\phi (A) \nsubseteq L$.\\
- $\forall a \in A, ~~ \phi(a) \equiv \gamma(a)\tau^0 \mod \tau$.\\
We denote most of the time $\phi_a$ instead of $\phi(a)$.
\item Given $\phi$ and $\psi$ two Drinfeld modules, an isogeny from $\phi$ to $\psi$ is a non-zero polynomial $f \in L\{\tau\}$ such that $f\cdot \phi_a = \psi_a \cdot f$ for all $a \in A$.
\end{itemize}
\end{df}
\begin{prop}[\cite{yu1995isogenies}]~
\begin{itemize}
\item Let $\phi$ and $\psi$ be Drinfeld modules over the $A$-field $L$. \\
If there exists an isogeny $f(\tau) \in L\{ \tau \}$ from $\phi \longrightarrow \psi$. Then there exists conversely also an isogeny $g \in L\{ \tau \}$ from $\psi \longrightarrow \phi$ such that\\
$f \cdot g = \psi_a \text{ and } g \cdot f = \phi_a \text{ for some } a \in A. $\\
This suggests that the isogeny relation is an equivalence relation.
\item let $s=[L : \mathbb{F}_q]$. $\pi_{\phi}=\tau^s : \phi \longrightarrow \phi $ is a special endomorphism  of $\phi$ called the Frobenius endomorphism. 
\item $\pi_{\phi}$ is an algebraic integer over $A$.
\item Two Drinfeld modules $\phi$ and $\psi$ over $L$ are isogenous iff the minimal polynomials of their Frobenius endomorphism $\pi_{\phi}$ and $\pi_{\psi}$ coincide i.e. $m_{\phi}=m_{\psi}$
\end{itemize}
\end{prop}
\begin{prop}[\cite{yu1995isogenies}]
Let $r=rank\phi$. The Frobenius enodmorphism $\pi_{\phi}$ satisfies the following:
\begin{enumerate}[1]
\item $\pi_{\phi}$ is integral over $A$ \label{c1}
\item There is only one place of the function field $k(\pi_{\phi})$ which is a zero of $\pi_{\phi}$ and this place lies above the place $v$.\label{c2}
\item There is only one place of $k(\pi_{\phi})$ lying over the place $\infty$ of $k$.\label{c3}
\item $\vert \pi_{\phi} \vert_{\infty}= l^{1/r}$ where $l=\vert L \vert$ \\
and $\vert . \vert_{\infty}$ is the unique extension to $k(\pi_{\phi})$ of the normalized absolute value of $k$ corresponding to the place $\infty$.\label{c4}
\item $[k(\pi_{\phi}): k]$ divides $r$\label{c5}
\end{enumerate}
\end{prop}
\begin{df}\label{Weil}
Any element $\pi \in \overline{k}$ satisfying the above properties is called \textbf{Weil number} and we call the corresponding minimal polynomial \textbf{Weil polynomial}. 
\end{df}
\begin{rmk}
We therefore completely describe the isogeny classes of rank $r$ Drinfeld modules by explicitly giving the list of all the corresponding Weil polynomials. 
\end{rmk}
\subsection{Separable Weil polynomials}
We assume here (as indicated by the title) that the function field extension $k(\pi)/k$ is separable. We come back later on to the inseparable case.
\begin{thm}\label{theo1} [characterization of Weil Polynomials]
Let $r \in \mathbb{N}$ (coprime to $q$). An irreducible polynomial $M(x) \in A[x]$ is  a ``rank" $r$ Weil polynomial if and only if it is of the form\\
\begin{equation}\label{poly}
M(x)=x^{r_1} +a_1x^{r_1-1} + \cdots + a_{r_1-1}x + \mu \mathfrak{p}_v^{\frac{m}{r_2}}
\end{equation}
such that the following conditions hold:
\begin{itemize}
\item[a)] $r_1, r_2 \in \mathbb{N}$ such that $r=r_1\cdot r_2$ and $r_2$ divides $m$. 
\item[b)] $\deg a_i \leq \frac{im\deg\mathfrak{p}_v}{r} $\\
\item[c)] $M_0(x)=x^{r_1} + \frac{a_1}{T^s}x^{r_1-1} + \cdots + \frac{a_{r_1-1}}{T^{s(r_1-1)}}x+ \mu\frac{\mathfrak{p}_v^{\frac{m}{r_2}}}{T^{sr_1}} \mod \frac{1}{T^h}$ is irreducible. Where $s,~h \in \mathbb{N}$ given by\\ 
$s=\lceil{\frac{m\deg\mathfrak{p}_v}{r}}\rceil$; $h=v_{\infty}\left(disc\left(M(x)\right)\right)+sr_1(r_1-1)+1~~$ \\
\item[d)] If $\overline{M(x)} \equiv \overline{f_1}(x) \cdot \overline{f_2}(x) \cdots \overline{f_{\mathfrak{s}}}(x) \mod \mathfrak{p}_v^n$ is an irreducible decomposition of $M(x) \mod \mathfrak{p}_v^n$ then $Res\left( \overline{f_i}(x), \frac{\overline{M(x)}}{\overline{f_i}(x)}\right) \not\equiv 0 \mod \mathfrak{p}_v~~~ \forall i=1, \cdots , \mathfrak{s}$.\\
 Where $n=v\left(disc\left(M(x)\right)\right)+1$.
\end{itemize}
%$r_1$ and $r_2$ here are positive integers such that $r=r_1\cdot r_2$ and $r_2$ divides $m$.
\end{thm}
\noindent We recall that we have assumed here $r$ to be coprime to $q$ so that the polynomial $M(x)$ remains separable.\\
Before proving the theorem, let us state the following lemmas, which will be useful for the proof.

\begin{lem}\label{prop6}
Let $k(\pi)$ be the function field generated by a root $\pi$ of an irreducible polynomial $M(x) \in k[x]$ of the form given in~\ref{poly}. Let  $M(x)=f_1(x) \cdot f_2(x) \cdots f_s(x)$ be the irreducible decomposition of $M(x)$ over $k_v$. If $f_{i_0}(x)$ describes a zero $\mathfrak{p}_{i_0}$ of $\pi$ in $k(\pi)$ (see \cite[Proposition 8.2]{Neukirch}), then $\pi$ has a unique zero in $k(\pi)$ if and only if $Res\left( f_{i_0}(x), \frac{M(x)}{f_{i_0}(x)} \right) \not\equiv~ 0 \mod \mathfrak{p}_v $.
$Res(?,?)$ denotes the resultant function.
\end{lem}
\begin{pf} 
Let us assume that $\pi$ has a unique zero $\mathfrak{p}_{i_0}$ in $k(\pi)$ described by the factor $f_{i_0}(x)$. If $Res\left(f_{i_0}(x), \frac{M(x)}{f_{i_0}(x)}\right) \equiv 0 \mod\mathfrak{p}_v$ then we have the following:\\
We recall that $\mathfrak{p}_v$ can be seen here as the unique place of the completion field $k_v$.\\
$\mathfrak{p}_v \mid Res\left(f_{i_0}(x), \frac{M(x)}{f_{i_0}(x)}\right)$ i.e. $\mathfrak{p}_v \mid Res\left(f_{i_0}(x), f_j(x)\right)$ for some $j \in \{1, \cdots, s  \}$\\$j\neq i_0$.  That means $\mathfrak{p}_i \mid Res\left(f_{i_0}(x), f_j(x)\right)$ for all $i=1, \cdots , s$.\\ 
Where $\mathfrak{p}_v = \mathfrak{p}_1^{n_1} \cdots \mathfrak{p}_s^{n_s}$ is the prime decomposition of $\mathfrak{p}_v$ in $k(\pi)$.\\
$\mathfrak{p}_i$ can be seen as the unique extension of $\mathfrak{p}_v$ in the completion field\\ $\left(k(\pi)\right)_{\mathfrak{p}_i}\simeq k_v(\pi_i)$. Where $\pi_i$ is a root of the irreducible factor $f_i(x) \in k_v[x]$ of $M(x)$ defining the place $\mathfrak{p}_i$.\\
In particular $\mathfrak{p}_{i_0}$ divides $Res\left(f_{i_0}(x), f_j(x)\right)$.\\
Let $\tilde{\mathfrak{p}}_{i_0}$ be a prime of $F$ above $\mathfrak{p}_{i_0}$. $F=Gal\left(k(\pi)\right)$ denotes the Galois closure of $k(\pi)$ (i.e. the splitting field of $M(x)$).\\
$\mathfrak{p}_{i_0} \mid Res\left(f_{i_0}(x), f_j(x)\right)$ implies that $\tilde{\mathfrak{p}}_{i_0} \mid Res\left(f_{i_0}(x), f_j(x)\right)$. In other words $\tilde{\mathfrak{p}}_{i_0}$ divides $\pi_{i_0}-\pi_j$ for some root $\pi_{i_0}$ of $f_{i_0}(x)$ and $\pi_j$ of $f_j(x)$.\\
$\tilde{\mathfrak{p}}_{i_0}$ divides $\pi$. $\pi$ and $\pi_{i_0}$ are both, roots of $f_{i_0}(x)$. The corresponding valuation $v_{i_0}$ is defined by $v_{i_0}=\overline{v}\circ \tau_0$ with
$$
\begin{tabular}{rcl}
$\tau_0: k(\pi)$ & $\xhookrightarrow{~~~~~}$ & $\overline{k_v}$\\
$\pi$ & $\longmapsto$ & $\pi_{i_0}$
\end{tabular}
$$
$\overline{v}$ is the valuation defined over $\overline{k_v}$ (extending $v$).\\
By definition, $\mathfrak{p}_{i_0}$ divides $\pi$ i.e. $v_{i_0}(\pi) > 0$. In addition, $\pi= \sigma(\pi_{i_0})$ for some $\sigma \in Gal(F/k)$. That is $v_{i_0}\circ \sigma (\pi_{i_0}) > 0$.\\
But $v_{i_0} \circ \sigma $ and $v_{i_0}$ define the same place of $k(\pi)$ because $\pi$ and $\pi_{i_0}$ are roots of the same irreducible factor $f_{i_0}(x)$. Thus $\tilde{\mathfrak{p}}_{i_0}$ divides $\pi_{i_0}.$ \\
$\tilde{\mathfrak{p}}_{i_0}$ divides $\pi_{i_0}-\pi_j$ and $\tilde{\mathfrak{p}}_{i_0}$ divides $\pi_{i_0}$ implies that $\tilde{\mathfrak{p}}_{i_0}$ divides $\pi_j$. \\
But $\pi_j=\sigma_j(\pi)$ for some $\sigma_j \in Gal(F/k)$. That means $\tilde{\mathfrak{p}}_{i_0} \mid \pi_j$ i.e. $\tilde{\mathfrak{p}}_{i_0} \mid \sigma_j(\pi)$. In other word $\sigma_j^{-1}\left(\tilde{\mathfrak{p}}_{i_0}\right) \mid \pi$.\\
$\sigma_j^{-1}\left(\tilde{\mathfrak{p}}_{i_0}\right)$ is a place of $F$ above the place $\mathfrak{p}_j$ of $k(\pi)$ defined by $f_j(x)$.\\
We have then $\pi \in \sigma_j^{-1}\left(\tilde{\mathfrak{p}}_{i_0}\right) \cap k(\pi) = \mathfrak{p}_j$.\\
Therefore $\pi$ possesses at least two zeros and it contradicts our initial hypothesis.\\
Let us assume conversely that $Res\left( f_{i_0}(x), \frac{M(x)}{f_{i_0}(x)} \right) \mod \mathfrak{p}_v \neq~ 0$.\\
If there are more than a zero of $\pi$ above $\mathfrak{p}_v$ in $k(\pi)$, then we have the following:\\
$M(x)= f_1(x) \cdots f_s(x) \in k_v[x]$. Suppose that $f_{i_0}(x)$ and $f_{i_1}(x)$ describe zeros of $\pi$ above $\mathfrak{p}_v$ in $k(\pi)$. Let $\tilde{\mathfrak{p}}_{i_0}$ and $\tilde{\mathfrak{p}}_{i_1}$ be primes of $F$ above $\mathfrak{p}_{i_0}$ and $\mathfrak{p}_{i_1}$ respectively. There exists $\sigma \in Gal(F/k)$ such that $\tilde{\mathfrak{p}}_{i_1}=\sigma\left(\tilde{
\mathfrak{p}}_{i_0}\right)$. Since $\tilde{\mathfrak{p}}_{i_0}$ and $\tilde{\mathfrak{p}}_{i_1}$ both divide $\pi$ we have $\sigma\left(\tilde{\mathfrak{p}}_{i_0}\right)$ divides $\pi$ and $\tilde{\mathfrak{p}}_{i_0}$ divides $\pi$. That is, $\tilde{\mathfrak{p}}_{i_0}$ divides $\sigma^{-1}(\pi)$ and $\tilde{\mathfrak{p}}_{i_0}$ divides $\pi$.\\
$\sigma^{-1}(\pi)$ is a conjugate of $\pi$ which is not a root of $f_{i_0}(x)$. Otherwise it would describe the same place of $k(\pi)$. Which is not the case since $\tilde{\mathfrak{p}}_{i_0}$ and $\tilde{\mathfrak{p}}_{i_1}$ are primes of $F$ above two distinct primes $\mathfrak{p}_{i_0}$ and $\mathfrak{p}_{i_1}$ of $k(\pi)$.\\
Thus $\tilde{\mathfrak{p}}_{i_0}$ divides $\sigma^{-1}(\pi) - \pi$. i.e. $\tilde{\mathfrak{p}}_{i_0}$ divides $Res\left(f_{i_0}(x), \frac{M(x)}{f_{i_0}(x)} \right)$. But $Res\left(f_{i_0}(x), \frac{M(x)}{f_{i_0}(x)} \right) \in A_v$. That is $Res\left(f_{i_0}(x), \frac{M(x)}{f_{i_0}(x)} \right) \in A_v \cap \tilde{\mathfrak{p}}_{i_0} = \mathfrak{p}_v$.\\
Therefore $\mathfrak{p}_v \mid Res\left(f_{i_0}(x), \frac{M(x)}{f_{i_0}(x)} \right)$ i.e. $Res\left(f_{i_0}(x), \frac{M(x)}{f_{i_0}(x)} \right) \equiv 0 \mod\mathfrak{p}_v$.\\
It contradicts our initial hypothesis.\\
Hence there is a unique zero of $\pi$ above $\mathfrak{p}_v$ in $k(\pi)$. 
\end{pf}

\begin{lem}\label{prop7}
Let $M(x)$ be as in the previous lemma. $\mathfrak{p}_1, \cdots , \mathfrak{p}_s$ denote the primes of $k(\pi)$ above $\mathfrak{p}_v$. If there is a unique prime containing $\pi$ i.e.
$$ \exists !~ i_0 \in \{1, \cdots , s\} \text{such that } \pi \in \mathfrak{p}_{i_0} \text{ but } \pi \notin \mathfrak{p}_j~ \forall j \neq i_0. $$
then so is it for any other conjugate $\tilde{\pi}$ of $\pi$.
\end{lem}

\begin{pf}
As mentioned before, $F$ denotes the splitting field of $M(x)$. $\pi$ and $\tilde{\pi}$ are conjugate. 
that means one can find $\alpha \in Gal\left( F/k \right)$ such that $\tilde{\pi}=\alpha\left( \pi \right)$.\\
$\pi \in \mathfrak{p}_{i_0}$ and $\pi \notin \mathfrak{p}_j ~\forall j \neq i_0$.
Let $\mathfrak{p}_{1j}, \cdots, \mathfrak{p}_{l_jj}$ be the primes of $F$ above $\mathfrak{p}_{j}$.\\ 
$\pi \in \mathfrak{p}_{i_0}$ means $\pi \in \mathfrak{p}_{ii_0} ~ \forall i=1, \cdots, l_{i_0}$. i.e. $\alpha(\pi) \in \alpha(\mathfrak{p}_{ii_0}) ~ \forall i=1, \cdots, l_{i_0}$. In other words $\tilde{\pi} \in \alpha(\mathfrak{p}_{i_0})$.\\
$\forall j \neq i_0 ~ \pi \not\in \mathfrak{p}_{j}$. That means $\pi \not\in \mathfrak{p}_{ij}$ for some $i \in \{1, \cdots, l_j\}$. Equivalently, $\alpha(\pi) \not\in \alpha(\mathfrak{p}_{ij})$ for some $i$. In other words $\tilde{\pi} \not\in \alpha(\mathfrak{p}_{j})$.\\
Therefore $\tilde{\pi} \in \alpha\left(\mathfrak{p}_{i_0} \right)$ and $\tilde{\pi} \notin \alpha\left(\mathfrak{p}_j \right) ~\forall j \neq i_0$. Since $\alpha$ acts as a permutation on the set of primes, we can conclude that $\tilde{\pi}$ belongs to some prime\\ $\mathfrak{q}_{k_0}=\alpha\left( \mathfrak{p}_{i_0}\right)$ of $k(\tilde{\pi})$ above $\mathfrak{p}_v$ and $\tilde{\pi}$ does not belong to any other prime $\mathfrak{q}_{j}~ j \neq k_0$ of $k(\tilde{\pi})$ above $\mathfrak{p}_v$.
\end{pf}
\begin{cor}\label{cor1}
$M(x)=f_1(x)\cdot f_2(x) \cdots f_s(x) \in k_v[x]$ is the irreducible decomposition in $k_v[x]$ of the polynomial $M(x)$. %of the first form~(\ref{p1}). \\
There is a unique zero of $\pi$ in $k(\pi)$ lying over the place $v$ of $k$ if and only if \\
$Res\left( f_{j}(x), \frac{M(x)}{f_{j}(x)} \right) \not\equiv 0 \mod \mathfrak{p}_v ~ \forall j \in \{1, \cdots, s\}$. 
%If $Res\left( f_{i_0}(x), \frac{M(x)}{f_{i_0}(x)} \right) \mod \mathfrak{p}_v \neq 0$ for some factor $f_{i_0}(x)$ of the above mentioned decomposition of $M(x)$ in $k_v[x]$ then\\
%$Res\left( f_{j}(x), \frac{M(x)}{f_{j}(x)} \right) \mod \mathfrak{p}_v \neq 0~ \forall j \in \{1, \cdots, s\}$.
\end{cor}
\begin{pf}
Let us assume that there is a unique zero of $\pi$ in $k(\pi)$. Let $f_{i_0}(x)$ be the irreducible factor of $M(x)$ in $k_v[x]$ describing that zero of $\pi$. That means
$Res\left( f_{i_0}(x), \frac{M(x)}{f_{i_0}(x)} \right) \mod \mathfrak{p}_v \neq 0$. If for some other $i_1 \in \{1, \cdots ,s\} ~ i_1 \neq i_0,~ Res\left( f_{i_1}(x), \frac{M(x)}{f_{i_1}(x)} \right) \mod \mathfrak{p}_v = 0$, then we have the following:\\
%$f_{i_0}(x)$ describes a zero $\mathfrak{p}_{i_0}$ in $k(\overline{\pi})$ of some root $\overline{\pi}$ of $M(x)$. We can therefore conclude from proposition~\ref{prop6} that $\mathfrak{p}_{i_0}$ is the unique zero of $\overline{\pi}$ in $k(\overline{\pi})$ above $\mathfrak{p}_v$ since $Res\left( f_{i_0}(x), \frac{M(x)}{f_{i_0}(x)} \right) \mod \mathfrak{p}_v \neq 0$. \\
$f_{i_1}(x)$ also describes a zero in $k(\tilde{\pi})$ of some root $\tilde{\pi}$ of $M(x)$. Let $F$ be the splitting field of $M(x)$. Since $M(x)$ is irreducible and separable over $k$, $Gal(F/k)$ acts transitively on the set of roots. That means  $\pi$ and $\tilde{\pi}$ are conjugate. In other words there exists $\alpha \in Gal(F/k)$ such that\\ 
$\tilde{\pi}=\alpha \left(\pi\right)$. $Res\left( f_{i_1}(x), \frac{M(x)}{f_{i_1}(x)} \right) \mod \mathfrak{p}_v=0$ means that $\tilde{\pi}$ has more than a zero in $k(\tilde{\pi} )$ above $\mathfrak{p}_v$ (see lemma~\ref{prop6}). This is (based on lemma~\ref{prop7}) a contradiction.\\
Hence we also have $Res\left( f_j(x), \frac{M(x)}{f_{j}(x)} \right) \mod \mathfrak{p}_v \neq 0 $ for any other $j \neq i_0.$\\
Conversely if $Res\left( f_j(x), \frac{M(x)}{f_{j}(x)} \right) \mod \mathfrak{p}_v \neq 0 $ for all $j \in \{1, \cdots , s\}$ then we have in particular $Res\left( f_{i_0}(x), \frac{M(x)}{f_{i_0}(x)} \right) \mod \mathfrak{p}_v \neq 0 $. Where $f_{i_0}(x)$ denotes an irreducible factor  of $M(x)$ in $k_v[x]$ describing a zero of $\pi$. Hence $\pi$ has a unique zero in $k(\pi)$ above the place $v$ of $k$ (see lemma~\ref{prop6}). 
\end{pf}
\vspace*{1cm}

\begin{pf}[of the theorem~\ref{theo1}]~~\\
\fbox{\parbox[b]{0.5cm}{$\Leftarrow$}} Under the assumption of our theorem, let $\pi$ be a rank $r$ Weil number and $M(x)$ be the correspoding Weil polynomial. We know from the condition~\ref{c5} of the definition of Weil number that $\deg M(x)=[k(\pi) : k]$ divides $r$. We denote then $r_1=\deg M(x)$ and $r_2=\frac{r}{r_1}$.\\
i.e. $M(x)$ has the form $M(x)=x^{r_1} +a_1x^{r_1-1} + \cdots + a_{r_1-1}x + a_{r_1}$.\\
Also $a_{r_1}=M(0)=(-1)^{r_1}N_{k(\pi)/k}\left( \pi \right)$. But $\pi$ has a unique zero in $k(\pi)$ which lies over $\mathfrak{p}_v$ according to the condition~\ref{c2} of the definition of Weil number. Thus $\mathfrak{p}_v$ is the unique prime of $A$ dividing $a_{r_1}$. That is 
\[ \hspace{5cm} a_{r_1}=\mu \mathfrak{p}_v^{\alpha} \hspace{5cm}(\star)  \label{eq*}\]
where $\alpha \in \mathbb{N}$, $\mu \in \mathbb{F}_q^*$\\
%and $R_v$ denotes the monic generator of the principal ideal $\mathfrak{p}_v$.\\
Moreover, we know from condition~\ref{c4} that $\vert \pi \vert_{\infty} = l^{1/r}=q^{\frac{m\deg \mathfrak{p}_v}{r}}$. That means $v_{\infty}(\pi)=v_{\infty}(\pi_i)=-\frac{m\cdot \deg\mathfrak{p}_v}{r}~~ \forall i ,\text{ where }\pi_{i's} \text{ denote the roots of } M(x)$. We also know that $a_{r_1}=\left(-1 \right)^{r_1}\displaystyle\prod_{i=1}^{r_1} \pi_i$. Hence $v_{\infty}(a_{r_1})=r_1v_{\infty}(\pi)=-\frac{m\cdot \deg\mathfrak{p}_v}{r_2}$. \\ 
 From \hyperref[eq*]{($\star$)} we have $-\alpha\deg\mathfrak{p}_v=v_{\infty}(a_{r_1} )=-\frac{m\cdot \deg\mathfrak{p}_v}{r_2}$ and therefore $\mathbb{N} \ni \alpha=\frac{m}{r_2}$. 
 $$ \text{Thus } r_2 \mid m \text{ and } a_{r_1}=\mu \mathfrak{p}_v^{\frac{m}{r_2}}$$ 
where $Q=\mathfrak{p}_v^m$ is the monic generator of the ideal~ $\mathfrak{p}_v^m$. Therefore $$M(x)=x^{r_1}+a_1x^{r_1-1} + \cdots + a_{r_1-1}x +\mu \mathfrak{p}_v^{\frac{m}{r_2}} \in A[x],~ \mu \in \mathbb{F}_q^*.$$
Let us consider again the roots $\pi_1, \cdots , \pi_{r_1}$ of $M(x)$ in $\overline{k}$. One knows that\\ $a_n=(-1)^n\displaystyle\sum_{i_1, \cdots , i_n}\pi_{i_1}\pi_{i_2} \cdots \pi_{i_n}$. That is\\ $v_{\infty}(a_n)=v_{\infty}\left( \displaystyle\sum_{i_1, \cdots , i_n}\pi_{i_1}\pi_{i_2} \cdots \pi_{i_n} \right) \geq \displaystyle\min_{i_1, \cdots , i_n} \Big\{ v_{\infty}(\pi_{i_1}\pi_{i_2} \cdots \pi_{i_n}) \Big\}$\\
But $$\displaystyle\min_{i_1, \cdots , i_n} \Big\{ v_{\infty}(\pi_{i_1}\pi_{i_2} \cdots \pi_{i_n}) \Big\}=v_{\infty}(\pi_{j_1}\pi_{j_2} \cdots \pi_{j_n})=v_{\infty}(\pi_{j_1}) + v_{\infty}(\pi_{j_2}) + \cdots +~v_{\infty}(\pi_{j_n})$$ $\text{ for some } (j_1, \cdots, j_n)$\\
Again as we mentioned before, one draws from condition~\ref{c4} that\\ $v_{\infty}(\pi_{j_1})=v_{\infty}(\pi_{j_2}) = \cdots =v_{\infty}(\pi_{j_n})=v_{\infty}(\pi)=-\frac{m\deg \frak{p}_v}{r}$. \\
Hence $v_{\infty}(a_n) \geq v_{\infty}(\pi_{j_1}) + v_{\infty}(\pi_{j_2}) + \cdots + v_{\infty}(\pi_{j_n})  = n\cdot v_{\infty}(\pi)= -\frac{n\cdot m\deg\mathfrak{p}_v}{r}$.\\
Thus $-\deg a_n \geq -\frac{n\cdot m\deg\mathfrak{p}_v}{r}$ that is $$\deg a_n \leq \frac{n\cdot m\deg\mathfrak{p}_v}{r}.$$
Therefore the coefficients $a_i$ of $M(x)$ satisfy the boundary condition $$\deg a_i \leq \frac{im\cdot\deg \mathfrak{p}_v}{r}=\frac{im\cdot\deg \mathfrak{p}_v^{\frac{m}{r_2}}}{r_1}.$$
Concerning the statement $c)$ of our theorem let us first of all recall that the irreducible polynomials\\
 $M(x)$ and $M_0(x)=x^{r_1} + \frac{a_1}{T^s}x^{r_1-1} + \cdots + \frac{a_{r_1-1}}{T^{s(r_1-1)}}x+ \mu\frac{\mathfrak{p}_v^{\frac{m}{r_2}}}{T^{sr_1}}$\\
define the same function field $k(\pi)=k\left(\frac{\pi}{T^s}\right)$.\\
One gets from \cite[Proposition 8.2]{Neukirch} that the decomposition of the place $\infty$ of $k$ in $k(\pi)$ is encoded in the decomposition of the polynomial $M_0(x)$ over the completion field $k_{\infty}$.\\
Therefore a unique place of $k(\pi)$ lying over the place at $\infty$ of $k$ if and only if the polynomial $M_0(x)$ is irreducible or a power of an irreducible polynomial over the completion field $k_{\infty}$. Since $r$ (and a fortiori $r_1$) is coprime to $q$, the polynomial $M_0(x)$ is separable. That means, there is a unique place of $k(\pi)$ over the place $\infty$ of $k$ if and only if $M_0(x)$ is irreducible over $k_{\infty}$.\\
In addition $h=v_{\infty}\left(disc\left(M(x)\right)\right)+sr_1(r_1-1)+1=v_{\infty}\left(disc\left(M_0(x)\right)\right)+1$.\\
One gets then from the Hensel lemma that any irreducible factor of $M_0(x) \mod \frac{1}{T^h}$ is the residue modulo $\frac{1}{T^h}$ of an irreducible factor of $M_0(x)$ in $k_{\infty}[x]$ and vice versa.\\
Hence there is a unique place of $k(\pi)$ lying over the place at $\infty$ of $k$ if and only if $M_0(x)$ is irreducible modulo $\frac{1}{T^h}$.\\
Concerning the statement $d)$ one has $n=v\left(disc\left(M(x)\right)\right)+1$.\\
 i.e. $disc\left(M(x)\right) \not\equiv 0 \mod \mathfrak{p}_v^n$.\\
One therefore concludes from corollary~\ref{cor1} in addition to the Hensel lemma that there is a unique (place) zero of $\pi$ in $k(\pi)$ if and only if\\
 $Res\left( \overline{f_i}(x), \frac{\overline{M(x)}}{\overline{f_i}(x)}\right) \not\equiv 0 \mod \mathfrak{p}_v~~~ \forall i=1, \cdots , \mathfrak{s}$.\\
 
 \fbox{\parbox[b]{0.5cm}{$\Leftarrow$}} Conversely let $\pi$ be a root of the polynomial in equation~\ref{poly} satisfying the four statements in theorem~\ref{theo1}. We aim to show that $\pi$ is a Weil number.\\
 First of all $\pi$ is an algebraic integer since $M(x)$ is a monic polynomial with coefficients in $A$.\\
 As we proved before, statement $c)$ implies that there is a unique place of $k(\pi)$ lying over the place at $\infty$ of $k$ and statement $d)$ implies that there is a unique (place) zero of $\pi$ in $k(\pi)$ and that place lies above the place $v$. \\
 Concerning the condition~\ref{c4} of the definition of Weil number, we have the following:\\
 from the constant coefficient (see equation~\ref{poly}) of the polynomial $M(x)$, one draws that the norm of $\pi$, $N_{k(\pi)/k}\left(\pi\right) = (-1)^{r_1}\mu\mathfrak{p}_v^{\frac{m}{r_2}} $. It therefore implies the following:\\
To avoid any ambiguity, let us denote $\infty '$ the place in $k(\pi)$ above $\infty$ in $k$. $v_{\infty '}\left( \pi \right):= \frac{1}{[k(\pi) : k]}v_{\infty}\left( N_{k(\pi)/k}\left(\pi \right) \right) = -\frac{1}{r_1}\deg \mathfrak{p}_v^{\frac{m}{r_2}} = -\frac{m\deg\mathfrak{p}_v}{r}$.\\
Therefore $\vert \pi \vert_{\infty '}=q^{-v_{\infty '}\left( \pi \right)}=q^{\frac{m\deg\mathfrak{p}_v}{r}}=l^{1/r}$.\\
The condition~ref{c5} follows from the statement $a)$ of our theorem~\ref{theo1}.\\
Hence $\pi$ is a Weil number and therefore its minimal polynomial $M(x)$ is a Weil polynomial.
  
\end{pf}

We summarize our result in the following algorithm, which one can use to check wether a given polynomial is a rank $r$ Weil polynomial or not.

\begin{algo}\label{algo}
\textbf{Input}: $M(x)=x^{r_1}+a_1x^{r_1-1} + \cdots + a_{r_1-1}x + \mu\mathfrak{p}_v^{\frac{m}{r_2}}$
\begin{enumerate}
\item Check that $r_1r_2=r$ and $\deg a_i \leq \frac{im\deg\mathfrak{p}_v}{r}$ for $i=1, 2, \cdots , r_1-1$. \\
If one of these conditions is not fulfilled then \textbf{Output} False and exit. \\ Else move to the next step.
\item Compute $D=disc\left(M(x)\right)$, $s=\lceil \frac{m\deg\mathfrak{p}_v}{r} \rceil$, $h= v_{\infty}\left( D \right) +sr_1(r_1-1) +1$ and $n=v(D)+1$. Where $v$ and $v_{\infty}$ are resp. the discrete valuations associated to the place $\mathfrak{p}_v$ and the place at $\infty$ of the field $k$. 
\item Set $M_0(x)=x^{r_1}+\frac{a_1}{T^s}x^{r_1-1}+ \frac{a_2}{T^{2s}} x^{r_1-2} + \cdots + \frac{a_{r_1-1}}{T^{s(r_1-1)}} x + \mu\frac{\mathfrak{p}_v^{\frac{m}{r_2}}}{T^{r_1.s}}$.\\
If $M_0(x)$ is not irreducible modulo $\frac{1}{T^h}$ then \textbf{output} False and exit.\\
else move to the next step. 
\item Compute $\overline{M}(x) \equiv M(x) \mod \mathfrak{p}_v^n$ and decompose (irreducibly)\\ $\overline{M}(x)=\bar{f}_1(x)\cdot \bar{f}_2(x) \cdots \bar{f}_{\mathfrak{s}}(x)$.\\ 
If for all $j \in \{1, \cdots , \mathfrak{s} \}$ $Res\left( \bar{f_j}(x), \frac{\overline{M}(x)}{\bar{f_j}(x)} \right) \neq 0 \mod \mathfrak{p}_v$ then\\ \textbf{output} True and exit.\\
Else: then \textbf{output} False and exit.\\
\end{enumerate}
\end{algo}
\begin{rmk}\label{rmk10}~
\begin{enumerate}
\item Each step of algorithm~\ref{algo} requires to know only the coefficients of the polynomial $M(x)$ and can be achieved in finitely many computations.
%\item As a corollary of the proposition~\ref{prop7}, it is actually not necessary to completely decompose $\overline{M}(x)$ modulo $\mathfrak{p}_v^n$. Only one irreducible factor is necessary for the test.
%\item If $r \mid \deg Q$ and $h=1$ then $s=\frac{\deg Q}{r}$, and the step 2 is done by simply checking that the polynomial\\ $M_0(x)=x^r + \displaystyle\sum_{i \in I} a_{i,0}x^{r-i} + \mu$ is irreducible over $\mathbb{F}_q$.\\Where $I=\Big\{ i=1, \cdots , r-1 ;~ \deg a_i = \frac{i\deg Q}{r} \Big\}$ and $a_{i,0}$ denotes the leading coefficient of $a_i$.\\
%Indeed, for $h=1$, $$ \frac{a_i}{T^{is}} \equiv \begin{cases}
%0 \mod \frac{1}{T} & if \deg a_i < \frac{i\deg Q}{r}\\
%a_{i,0} \mod \frac{1}{T} & if \deg a_i = \frac{i\deg Q}{r}
%\end{cases}$$
%One can also remark that the residue field associated to the place $\infty$ is $\mathbb{F}_q$.\\ One may also notice that $h=1$ if and only if $ disc \left( M_0(x) \right) \neq 0$ where $M_0(x)=x^r + \displaystyle\sum_{i \in I} a_{i,0}x^{r-i} + \mu$.\\
%Indeed, The discriminant of the polynomial $$M(x)=x^r +a_1x^{r-1} + \cdots + a_{r-1}x +\mu Q$$ is a homogeneous polynomial of degree $2r-2$ in its coefficients. One of the monomials is $\mu^{r-1} Q^{r-1}$ whose degree $(\text{in } T)$ is $(r-1)\deg Q$. Also, all the monomials of the form $a_{i_1}^{\alpha_1} \cdots a_{i_l}^{\alpha_l}$ with $i_k \in I,$ have degree $(\text{in } T)$\\ $(r-1)\deg Q$. Therefore $v_{\infty}\left( disc\left( M(x) \right)\right):=-\deg_T\left(disc\left(M(x)\right)\right)=-(r-1)\deg Q$ iff $disc\left( M_0(x) \right) \neq 0$. In such a case $h=1$.

\item A priori, the algorithm only decides for a given polynomial whether it is a Weil polynomial or not. But one can also use that algorithm to provide the complete list of rank $r$ Weil polynomials.\\
Indeed, the coefficients $a_i$ of the potential Weil polynomial $$M(x)=x^{r_1} + a_1x^{r_1-1} + \cdots + a_{r_1-1}x + \mu\mathfrak{p}_v^{\frac{m}{r_2}}$$ are bounded by $\deg a_i \leq \frac{im\deg\mathfrak{p}_v}{r}$ and $a_i \in \mathbb{F}_q[T].$ So there are finitely many such polynomials. One can then check for each such polynomial (using the algorithm) whether it is a rank $r$ Weil polynomial or not.\\
In fact the number of polynomials $a_i \in \mathbb{F}_q[T]$ of degree atmost $\frac{im\deg\mathfrak{p}_v}{r}$ is $q^{\frac{im\deg \mathfrak{p}_v}{r}+1}$. Thus for polynomials of the form $$M(x)=x^{r_1} + a_1x^{r_1-1} + \cdots + a_{r_1-1}x + \mu \mathfrak{p}_v^{\frac{m}{r_2}} \in A[x],$$ we have a total number of $\displaystyle\prod_{i=1}^{r_1-1} q^{\frac{im\deg\mathfrak{p}_v }{r}+1}=q^{(r_1-1)\left[1+\frac{m\deg\mathfrak{p}_v}{2r_2}\right]}$ polynomials to be checked. This number can be reduced if one takes into account the following result.
%%%%%%%%%%%%To be continued%%%%%%%%%%%%%%%%%%%%%%%%
%%%%%%%%%%%%%%%%%%%%%%%%%%%%%%%%%%%%%%%%%%%%%%%%%%%
\end{enumerate}
\end{rmk}
 \begin{prop}\label{c2-cor}
We consider the same polynomial\\ $M(x)=x^{r_1} + a_1x^{r_1-1} + \cdots + a_{r_1-1}x + \mu\mathfrak{p}_v^{\frac{m}{r_2}} \in A[x],$ whose root $\pi$ generates the function fields extension $k(\pi)/k$. If $\mathfrak{p}_v$ does not divide the linear coefficient $a_{r_1-1}$, then $\pi$ satisfies the condition~\ref{c2} of the definition of a Weil number. That is, there is a unique zero of $\pi$ in $k(\pi)$ over the place $v$.
\end{prop}

\begin{pf}
We proceed by contraposition of the above statement. That is, if $\pi$ has more than a zero over the place $v$ then $\mathfrak{p}_v$ divides $a_{r_1-1}$.\\
Let $\mathfrak{p}_1$ be a zero of $\pi$ above $v$ in $k(\pi)$. If $\pi$ has another zero say $\mathfrak{p}_2$, then we have the following.\\
Let $F$ be the splitting field of $M(x)$. $F/k$ is a Galois extension and $k(\pi)$ is an intermediate field. Let $\mathfrak{p}_1'$ and $\mathfrak{p}_2'$ be extensions of $\mathfrak{p}_1$ and $\mathfrak{p}_2$ respectively in $F$. Let $B$ be the integral closure of $A$ in $k(\pi)$. $\mathfrak{p}_1'\cap A=\mathfrak{p}_1'\cap B\cap A=\mathfrak{p}_1 \cap A=\mathfrak{p}_v$. Same for $\mathfrak{p}_2'$. So $\mathfrak{p}_1'$ and $\mathfrak{p}_2'$ are primes of $F$ above $\mathfrak{p}_v$. Since $Gal(F/k)$ acts transitively on the sets of primes above $\mathfrak{p}_v$, there exists $\sigma \in Gal(F/k)$ such that $\mathfrak{p}_2'=\sigma(\mathfrak{p}_1')$. $\mathfrak{p}_2' \big\vert \pi$ then $\sigma(\mathfrak{p}_1') \big\vert \pi $. That is $\mathfrak{p}_1' \big\vert \sigma^{-1} (\pi)$. Moreover, $\sigma^{-1}(\pi) \neq \pi$ otherwise $\sigma$ would be in $Gal(F/k(\pi))$ that is\\ $\mathfrak{p}_1=\sigma(\mathfrak{p}_1)=\sigma(\mathfrak{p}_1'\cap B)=\sigma(\mathfrak{p}_1')\cap \sigma(B)=\mathfrak{p}_2'\cap B=\mathfrak{p}_2$. Which is not possible since $\mathfrak{p}_1 \neq \mathfrak{p}_2$. Taking into account the following representation of the linear coefficient in terms of the roots of $M(x)$, that is $a_{r_1-1}=\displaystyle\sum_{j=1}^{r_1}\prod_{i=1,i\neq j}^{r_1} \tau_i(\pi)$ and in addition to the fact that $\mathfrak{p}_1' \big\vert \pi$ and $\mathfrak{p}_1' \big\vert \sigma^{-1}(\pi)$, we therefore get $\mathfrak{p}_1' \big\vert a_{r_1-1}.$ That is $a_{r-1} \in \mathfrak{p}_1'$ but $a_{r-1} \in A$. Hence $a_{r-1} \in \mathfrak{p}_1'\cap A=\mathfrak{p}_v$ i.e. $\mathfrak{p}_v \big\vert a_{r_1-1}$. \hspace{12cm}
\end{pf}
\begin{rmk}\label{rmk11}
As we mentioned in remark~\ref{rmk10}, if one takes into account the above mentioned result, the number of polynomials to be checked (using the whole algorithm~\ref{algo}) can be reduced to $$q^{1+\frac{(r_1-1)m\deg\mathfrak{p}_v}{r}-\deg \mathfrak{p}_v}\times \displaystyle\prod_{i=1}^{r_1-2} q^{1+\frac{im\deg\mathfrak{p}_v}{r}}=q^{(r_1-1)\left[ \frac{m\deg\mathfrak{p}_v}{2r_2}+1 \right]- \deg \mathfrak{p}_v}$$ 
For other polynomials for which $\mathfrak{p}_v \nmid a_{r_1-1}$, one can just check the step 2 of our algorithm.
\end{rmk}
\begin{rmk}
One easily notices that if $m$ and $r$ are coprime, then the only potential rank $r$ Weil polynomials are the one of the form
\begin{center}
$M(x)= x^{r} + a_1 x^{r-1} + \cdots + a_{r-1}x +\mu\mathfrak{p}_v^m$ 
\end{center}
%such that the polynomial $M_0(x)=x^{r} + \frac{a_1}{T^s}x^{r-1} + \cdots + \frac{a_{r-1}}{T^{s(r-1)}}x+ \mu\frac{Q}{T^{sr}}$ is irreducible in $k_{\infty}[x]$.
\end{rmk}
\subsection{Inseparable Weil polynomials}
As we mentioned before, the condition $r$ coprime to $q$ was made so that the polynomial $M(x)$ remains separable. Let us now drop that conition and pick any positive integer $r$.
\begin{rmk}
Before going further, let us draw the attention of the reader on the following fact:\\
The only sprain to the generality is how to check the conditions ~\ref{c2}~and~\ref{c3} when $M(x)$ is inseparable. In other words how to get the irreducible factorization of $M(x)$ over the completion field $k_* \in \{ k_{\infty},~ k_v \}$. In the previous case, the factorization was entirely determine by the irreducible decomposition of $M(x) \mod \mathfrak{p}_v^n$ and $M(x) \mod \frac{1}{T^h}$ for $k_v$ and $k_{\infty}$ respectively. Where $n=v\left( disc\left( M(x) \right)\right) + 1$ and $h=v_{\infty}\left( disc\left( M(x) \right) \right) +sr_1(r_1-1) +1.$\\
That argument is not valid anymore in this case because $disc \left( M(x) \right) = 0$. But at least one knows that if $M(x)$ is an inseparable irreducible polynomial over a field $k$ of characteristic $p > 0$, then there exists a separable polynomial $f(x) \in k[x]$ such that $M(x)=f\left( x^{p^d} \right)$ for some $d \in \mathbb{N}$. We will use the separable polynomial $f(x)$ to overcome the difficulties encountered when $M(x)$ is inseparable.  
\end{rmk}
\subsubsection*{Some properties of monic irreducible polynomials over a field k of characteristic $p > 0$}
We provide in this part, as mentioned in the title, some important properties of irreducible polynomials over a field $k$, with $char (k)=p >0$. These properties will be very helpful later on.

\begin{prop}\cite[theorem A6, page 11]{KConrad}\label{prop0}~\\
Let $k$ be a field of characteristic $p >0$ and $f(x)$ be a monic irreducible polynomial in $k[x]$. Then $f(x^p)$ is either irreducible or a $p$-th power of an irreducible polynomial in $k[x]$. 
\end{prop} 
\begin{pf}\cite{KConrad} 
\end{pf}
\begin{cor}\cite[Corollary A8]{KConrad}\label{cor0}~\\
Let $k$ be a field of characteristic $p > 0$ and $f(x)$ be a monic irreducible polynomial in $k[x]$. The following statements are equivalent.
\begin{enumerate}
\item[(i)] $f(x^{p^n})$ is irreducible in $k[x] ~ \forall n \in \mathbb{N}$.
\item[(ii)] $f(x) \notin k^p[x]$
\end{enumerate}
One should keep in mind that we mean by $k^p=\{ a^p, ~ a \in k\}$.
\end{cor}
\begin{pf}\cite{KConrad}~
\end{pf}
\begin{cor}\label{cor00}
Let $k$ be a field of characteristic $p >0$ and $f(x)$ be a monic irreducible polynomial in $k[x]$. Let $n \in \mathbb{N}$. \\
$f(x^{p^n})$ is either irreducible or a $p^{n_0}$-th power of an irreducible polynomial in $k[x]$ for some $n_0 \in \mathbb{N}$.
\end{cor}
\begin{pf}
Let $f(x)$ be a monic irreducible polynomial in $k[x]$ as mentioned in the corollary above. We know from corollary~\ref{cor0} that if $f(x) \notin k^p[x]$ then $f(x^{p^n})$ is irreducible.\\
Now if $f(x) \in k^p[x]$ then, \\
Let $f(x)=x^d + a_1^px^{d-1}+ \cdots + a_{d-1}^px + a_d^p$.\\
We set $n_0=\min \Big\{ \nu_p(a_i^p),~i=1, \cdots , d \Big\}$ where $\nu_p(a_i^p)$ denotes the positive integer $t$ such that $a_i^p= b_i^{p^t}$ and $b_i \in k \setminus{k^p}$. Let $a_{i_0}^p$ be the coefficient for which $n_0=\nu_p(a_{i_0}^p).$\\
$$f(x)=x^d + b_1^{p^{n_0+r_1}}x^{d-1}+ b_2^{p^{n_0+r_2}}x^{d-2} + \cdots + b_{i_0}^{p^{n_0}}x^{d-i_0}+ \cdots + b_{d-1}^{p^{n_0+r_{d-1}}}x + b_d^{p^{n_0+r_d}}$$
\begin{itemize}
\item[$If ~$]$n \geq n_0$ then we have the following\\
$
\begin{array}{rcl}
f(x^{p^n}) & = & x^{d p^n} + b_1^{p^{n_0+r_1}}x^{(d-1)p^n}+\cdots + b_{i_0}^{p^{n_0}}x^{(d-i_0)p^n} + \cdots +\\
& & + b_{d-1}^{p^{n_0+r_{d-1}}}x^{p^n} + b_d^{p^{n_0+r_d}}\\
 & = & \left( x^{dp^{n-n_0}} + b_1^{p^{r_1}}x^{(d-1)p^{n-n_0}}+ \cdots + b_{i_0}x^{(d-i_0)p^{n-n_0}} + \cdots +\right. \\
 & & + \left. b_{d-1}^{p^{r_{d-1}}}x^{p^{n-n_0}} + b_d^{p^{r_d}}\right)^{p^{n_0}} \\
 & = & \left( g_0\left(x^{p^{n-n_0}} \right) \right)^{p^{n_0}}
\end{array}
$\\
with $g_0(x)=x^{d} + b_1^{p^{r_1}}x^{d-1}+ \cdots + b_{i_0}x^{d-i_0} + \cdots + b_{d-1}^{p^{r_{d-1}}}x + b_d^{p^{r_d}}$\\
$g_0(x)$ must be irreducible in $k[x]$. Indeed,\\
If $g_0(x)$ is reducible in $k[x]$, that is $g_0(x)=h_1(x)\cdot h_2(x)$ with $h_1(x) \text{ and } h_2(x)$ in $k[x]$, then we have the following:\\
$g_0\left(x^{p^{n-n_0}} \right)=h_1\left(x^{p^{n-n_0}} \right) \cdot h_2\left(x^{p^{n-n_0}} \right)$. That is, \\ 
$
\begin{array}{rcl}
f(x^{p^n})=\left(g_0\left(x^{p^{n-n_0}} \right) \right)^{p^{n_0}} & = &\left(h_1\left(x^{p^{n-n_0}} \right)\right)^{p^{n_0}} \cdot \left(h_2\left(x^{p^{n-n_0}} \right)\right)^{p^{n_0}}\\
& = & h_1^{p^{n_0}}\left( x^{p^n}\right) \cdot h_2^{p^{n_0}}\left( x^{p^n}\right)
\end{array}
$\\
Where $ h_i^{p^{n_0}}\left(x\right)$ denotes the polynomial obtained from $h_i(x)$ by raising all its coefficients to the power $p^{n_0}$.\\
Thus $f(x^{p^n})=h_1^{p^{n_0}}\left( x^{p^n}\right) \cdot h_2^{p^{n_0}}\left( x^{p^n}\right)$ i.e. $f(x)=h_1^{p^{n_0}}\left( x\right) \cdot h_2^{p^{n_0}}\left( x \right)$ which contradicts the fact that $f(x)$ is irreducible.\\
Hence $g_0(x)$ must be irreducible in $k[x]$.\\
In addition, since $b_{i_0} \notin k^p$, we also have $g_0\left( x^{p^{n-n_0}} \right)$ is irreducible (see corollary~\ref{cor0}).
\item[$If $]$n < n_0$ then one can write down $f(x^{p^n})$ as follows $$F(x):=f(x^{p^n})=~\left( g(x) \right)^{p^n}$$ 
with $g(x)=x^d + c_{d-1}x^{d-1} + \cdots + c_1 x + c_0 \in k^p[x]$.\\
\underline{Claim 1}: If $f(x)$ is separable then so is $g(x)$.\\
We know that $f(x)$ is a separable polynomial and $d= \deg f(x)$. We also know that for each root $\alpha$ of $f(x)$, the $p^n$-th root $\alpha^{\frac{1}{p^n}}$ of $\alpha$ is a root of $F(x)$. So $F(x)$ has at least $d$ distinct roots \circled{1}.\\
Also $F(x)=f(x^{p^n})=\left( g(x) \right)^{p^n}$ and $\deg g(x)=d$. Thus $F(x)$ has a maximum of $d$ distinct roots \circled{2}.\\
\circled{1} and \circled{2} imply that $F(x)$ must have exactly $d$ distinct roots.\\
Therefore $g(x)$ is separable.\\
\underline{Claim 2}: $g(x)$ is irreducible over $k$.\\
Indeed, Let us assume that $g(x)$ is reducible over $k$.\\
That is $g(x)=h_1(x) \cdot h_2(x)$.\\
Therefore $\left( g(x) \right)^{p^n}=\left(h_1(x)\right)^{p^n} \cdot \left(h_2(x)\right)^{p^n}=h_1^{p^n}\left( x^{p^n} \right) \cdot h_2^{p^n}\left( x^{p^n} \right)$.
Where $h_i^{p^n}\left( x \right)$ denotes the polynomial obtained from $h_i(x)$ by raising all its coefficients to the power $p^n$.\\
Thus $f(x^{p^n})=\left( g(x) \right)^{p^n}=h_1^{p^n}\left( x^{p^n} \right) \cdot h_2^{p^n}\left( x^{p^n} \right)$ That is $f(x)=h_1^{p^n}\left( x \right) \cdot~ h_2^{p^n}\left( x \right)$
which is impossible since $f(x)$ is irreducible over $k$.\\
Hence $g(x)$ must be irreducible.\\
So for this special case, if in addition to the hypothesis of the corollary $f(x)$ is separable, then $f(x^{p^n})$ would be a $p^n$-th power of an irreducible separable polynomial.
\end{itemize}
Therefore in any case $f(x^{p^n})$ is either irreducible or a $p^{n_0}$-th power of an irreducible polynomial in $k[x]$. \hspace{8.5cm}
\end{pf}
~\\

\noindent Let us come back to our Weil number $\pi$ with all the notations we have set at the beginning and $k=\mathbb{F}_q(T)$. We now assume that the extension $k(\pi)/k$ is  not separable. That is the minimal polynomial $M(x)$ of $\pi$ is an irreducible inseparable polynomial in $k[x]$. We know that if it is the case, then there exists a separable irreducible polynomial $f(x) \in k[x]$ such that $$M(x) =~f(x^{p^n}) \text{ for some } n \in \mathbb{N}.$$
Let us first discuss the case where $n=1$. i.e. $M(x)=f(x^p)$.\\
Let $f(x)=f_1(x) \cdots f_{\mathfrak{s}}(x)$ be the irreducible decomposition of $f(x)$ over the completion field $k_*$ (where $k_* \in \big\{ k_v,~k_{\infty} \big\}$).\\
So $M(x)=f(x^p)=f_1(x^p) \cdots f_{\mathfrak{s}}(x^p)$. According to the proposition~\ref{prop0}, each polynomial $f_i(x^p)$ is either irreducible or a $p$-th power of an irreducible polynomial $h_i(x) \in k_*[x]$ i.e. $f_i(x^p) =\left( h_i(x)\right)^p$. In any case, the irreducible decomposition of $f(x)$ encodes all the irreducible factors of $M(x)$ in $k_*[x]$ and is enough to decide about the conditions~\ref{c2}~and~\ref{c3} of definition~\ref{Weil} of Weil number. Indeed, \\
$\pi$ satisfies condition~\ref{c3} if and only if $M(x)$ is irreducible or a power of an irreducible polynomial over $k_{\infty}$.\\
But we can say from our above discussion that $M(x)$ is irreducible or a power of an irreducible polynomial over $k_{\infty}$ if and only if the separable polynomial $f(x)$ is irreducible over $k_{\infty}$\\
%Likewise, $\pi$ satisfies condition $(\hyperref[c2]{c2})$ of definition~\ref{def1} if and only if there exists a unique irreducible factor $f_{i_0}(x)$ of $f(x)$ over the completion field $k_v$ such that $f_{i_0}(x^p)$ describes a zero $\mathfrak{p}_{i_0}$ of $\pi$ in $k(\pi)$. In other words such that $\mathfrak{p}_{i_0} \mid b_{0,i_0}$ and $Res\left(M(x), \frac{M(x)}{f_{i_0}(x^p)} \right) \neq 0 \mod \mathfrak{p}_v$ (see proposition~\ref{prop7} ). Where $b_{0,i_0}$ denotes the constant coefficient of the polynomial $f_{i_0}(x)$.\\
Likewise, one can properly check in this case the condition~\ref{c2}) of the definition of Weil number using proposition~\ref{prop6}, where $M(x)$ is replaced by the irreducible separable polynomial $f(x)$. In other words the condition~\ref{c2} is satisfied by the polynomial $M(x)$ if and only if it is satisfied by the polynomial $f(x)$. That is 
%$Res\left(f_1(x), \frac{f(x)}{f_1(x)} \right) \neq 0 \mod \mathfrak{p}_v$ which is equivalent to saying that
$Res\left(f_i(x), \frac{f(x)}{f_i(x)} \right) \neq 0 \mod \mathfrak{p}_v ~ \forall i \in \{1, \cdots , \mathfrak{s} \}$. Thanks to corollary~\ref{cor1}.\\
    
Now if $M(x)=f(x^{p^n})$ with $n > 1$ then the same idea holds. That~is, the irreducible decomposition of $f(x)=f_1(x) \cdots f_s(x)$ over the completion field $k_*$ encodes the irreducible decomposition of $M(x)$ over $k_*$. $$M(x)=f(x^{p^n})=f_1(x^{p^n}) \cdots f_{\mathfrak{s}}(x^{p^n})$$ 
From corollary~\ref{cor00}, one can draw that each $f_i(x^{p^n})$ is either irreducible or a $p^{n_0}$-th power of an irreducible polynomial in $k_*[x]$ for some $n_0 \in \mathbb{N}$.\\
Therefore one can use the irreducible decomposition of $f(x)$ in $k_*[x]$ to check the conditions~\ref{c2}~and~\ref{c3} of definition~\ref{Weil}. Exactly as it happened for the case $n=1$, \\
$\pi$ satisfies condition~\ref{c3} if and only if $M(x)$ is irreducible or a power of an irreducible polynomial over $k_{\infty}$.\\
$M(x)=f(x^{p^n})$ is irreducible or a power of an irreducible polynomial over $k_{\infty}$ if and only if the separable polynomial $f(x)$ is irreducible over $k_{\infty}$. Thanks to corollary~\ref{cor00}.\\
Following the same idea, the polynomial (or a root $\pi$ of the polynomial) $M(x)$ satisfies the condition~\ref{c2}) of definition~\ref{Weil} if and only if\\
$Res\left(f_i(x^{p^n}), \frac{M(x)}{f_i(x^{p^n})} \right) \neq 0 \mod \mathfrak{p}_v ~ \forall i \in \{1, \cdots , \mathfrak{s} \}$. Thanks once more to proposition~\ref{prop6} and also to corollary~\ref{cor1}. %which assure us that it is equivalent to saying that 
%$Res\left(f_1(x), \frac{f(x)}{f_1(x)} \right) \neq 0 \mod \mathfrak{p}_v$. In other words $Res\left(f_1(x^{p^n}), \frac{M(x)}{f_1(x^{p^n})} \right) \neq 0 \mod \mathfrak{p}_v$. Thanks once more to proposition~\ref{prop6} and also to corollary~\ref{cor1} which assure us that it is equivalent to saying that  
%Following the same idea, $\pi$ satisfies condition $(\hyperref[c2]{c2})$ of definition~\ref{def1} if and only if there exists a unique irreducible factor $f_{i_0}(x)$ of $f(x)$ over the completion field $k_v$ such that $f_{i_0}(x^{p^n})$ describes a zero $\mathfrak{p}_{i_0}$ of $\pi$ in $k(\pi)$. In other words such that $\mathfrak{p}_{i_0} \mid b_{0,i_0}$ and $Res\left(M(x), \frac{M(x)}{f_{i_0}(x^{p^n})} \right) \neq 0 \mod \mathfrak{p}_v$ (see proposition~\ref{prop7} ). Where $b_{0,i_0}$ denotes the constant coefficient of the polynomial $f_{i_0}(x)$. Thanks once more to corollary~\ref{cor00}. One can also (as we did for the separable case) mention corollary~\ref{cor_prop6} (using the separable polynomial $f(x)$ ) to justify the fact that the latter statement does not depend on the chosen root $\pi$ of $M(x)$. We mean that given two distinct conjugates roots $\pi$ and $\tilde{\pi}$ of $M(x)$, $\pi$ has a unique zero over $\mathfrak{p}_v$ in $k(\pi)$ if and only if $\tilde{\pi}$ has a unique zero over $\mathfrak{p}_v$ in $k(\tilde{\pi})$. 
\begin{rmk}
A conclusion one can draw from our discussion above is that, modulo some slight changes, one can use the same algorithm~\ref{algo} in the case where the polynomial $M(x)$ is inseparable. After those minor changes, we get the following algorithm.
\end{rmk} 
\begin{algo}\label{algo2}~ \\
$\mathbf{Input}: ~M(x)=x^{p^nr_0} + a_1x^{p^n(r_0-1)} + \cdots + a_{r_0-1}x^{p^n} + \mu \mathfrak{p}_v^{\frac{m}{r_2}}=f(x^{p^n})$.\\
Where $f(x)=x^{r_0} + a_1 x^{r_0-1} + \cdots + a_{r_0-1}x + \mu\mathfrak{p}_v^{\frac{m}{r_2}}$ with $r_0$ coprime to $q$.
\begin{enumerate}
\item Check that $p^nr_0r_2=r$ and $\deg a_i \leq \frac{im\deg\mathfrak{p}_v}{r}$ for $i=1, 2, \cdots, r_0-1$.\\
If one of these conditions does not hold then \textbf{output:} False and exit.\\
Else move to the next step.
\item Compute $D=disc\left(f(x)\right)$, $s=\ceil*{\frac{m\deg\mathfrak{p}_v}{r}}=\ceil*{\frac{m\deg\mathfrak{p}_v}{p^nr_0r_2}}$,\\
$h= v_{\infty}\left( D \right) +sr_0(r_0-1) +1$ and $u=v\left(disc\left(f(x)\right)\right)+1$, where $v$ and $v_{\infty}$ denote the discrete valuation associated to the place $\mathfrak{p}_v$ and the place at $\infty$ of the field $k$ respectively.
\item Set $f_0(x)=x^{r_0}+\frac{a_1}{T^{p^ns}}x^{r_0-1}+ \frac{a_2}{T^{2p^ns}} x^{r_0-2} + \cdots + \frac{a_{r_0-1}}{T^{p^n(r_0-1)s}} x + \mu\frac{\mathfrak{p}_v^{\frac{m}{r_2}}}{T^{p^nr_0.s}}$.\\
If $f_0(x)$ is not irreducible modulo $\frac{1}{T^h}$ then \textbf{output} False and exit.\\
else move to the next step.
\item Compute $\overline{f}(x) \equiv f(x) \mod \mathfrak{p}_v^u$ and provide the irreducible decomposition $\overline{f}(x)=\bar{f}_1(x)\cdot \bar{f}_2(x) \cdots \bar{f}_{\mathfrak{s}}(x)$. That is the irreducible decomposition of $\overline{M}(x)$ is given by $\overline{M}(x)=\overline{f}(x^{p^n})=\bar{f}_1(x^{p^n})\cdot \bar{f}_2(x^{p^n}) \cdots \bar{f}_{\mathfrak{s}}(x^{p^n})$\\
If for all $j \in \{1, \cdots \mathfrak{s} \}$ $Res\left(\bar{f_j}(x), \frac{\bar{f}(x)}{\bar{f_j}(x)} \right) \neq 0 \mod \mathfrak{p}_v$ \\
%i.e.\\ $Res\left(\bar{f_j}(x^{p^n}), \frac{M(x)}{\bar{f_j}(x^{p^n})} \right) \neq 0 \mod \mathfrak{p}_v$
then \textbf{output} True and exit.\\
Else: \textbf{output} False and exit. 
\end{enumerate} 
%\item Let $b_{0,i}$ be the constant coefficient of $\bar{f_i}(x^{p^n})$ which is the same as the one of $\bar{f_i}(x)$.\\
%$i=1$\\ While $\mathfrak{p}_v \nmid b_{0,i}$ or $ Res\left( \bar{f_i}(x^{p^n}), \frac{\overline{M}(x)}{\bar{f_i}(x^{p^n})} \right) = 0 \mod \mathfrak{p}_v$; $~~i=i+1$.\\
%If $i==\mathfrak{s}+1$ then \textbf{Outputs} NO.\\
%else: \textbf{Outputs} YES 
%\\ $Res(\cdot, \cdot)$ denotes the resultant function.
\end{algo}
\begin{rmk}
The above mentioned algorithm is based on the fact that the irreducible decomposition of the separable polynomial $f(x)$ in $k_*[x]$ encodes the irreducible decomposition of $M(x)=f(x^{p^n})$ in $k_*[x]$. We mean that one can get a 1-to-1 map between the irreducible factors of $f(x)$ and those of $M(x)=f(x^{p^n})$ in $k_*[x]$. 
\end{rmk}
\section{Isomorphism classes}
We keep the same data we had before and we consider Drinfeld $A$-modules defined over the finite $A$-field $L$.\\ 
It is known in the Drinfeld modules theory that two Drinfeld modules are isomorphic over an algebraic closure $\overline{L}$ of $L$ if and only if the have the same $J$-Invariants. In other words the $J$-Invariants determine the $\overline{L}$-Isomorphism classes of Drinfeld modules defined over $L$.\\
A natural question to ask is: how about $L$-Isomorphism classes? How does one check that two Drinfeld modules are isomorphic over the field $L$ itself?\\
We give an answer to that question in this part by coming up with an additional invariant we call \textbf{Fine Isomorphy Invariant} and we prove that, that invariant together with the $J$-Invariants determine the $L$-Isomorphism classes of Drinfeld modules over $L$.

\begin{df}[Fine Isomorphy Invariant]~\\
\indent Let $\phi: A \longrightarrow L\{\tau\}$ be a rank r Drinfeld $A$-module defined by $$ \phi_T= \gamma(T)+ g_1\tau + \cdots + g_r\tau^r $$ We set \\
$$d=gcd(q^k-1, ~ k\in I)=q^{\delta}-1$$
$ \text{ where } I= \{ i=1, \cdots, r; ~g_i \neq 0\}
 \text{ and } \delta=gcd(k:~k\in I)$.\\
We write $d=\displaystyle\sum_{k \in I} \lambda_k(q^k-1); ~ \lambda_k \in \mathbb{Z}$ and we set $\lambda=(\lambda_k)_{k \in I}$.\\
Let $B=\left\lbrace\alpha=(\alpha_k)_{k \in I},~~d=\displaystyle\sum_{k \in I} \alpha_k(q^k-1) \right\rbrace$.\\ The fine isomorphy invariant of $\phi$ is defined as $FI(\phi)=\left(FI_{\lambda}(\phi) \right)_{\lambda \in B}$, where
 $$ FI_{\lambda}(\phi)= \displaystyle \prod_{k\in I} g_k^{\lambda_k} ~mod L^{*d}$$
\end{df}
\begin{ex}
Let $\phi: A \longrightarrow L\{\tau \}$ be a rank 2 Drinfeld module defined by $\phi_T=\gamma(T)+g_1\tau +g_2\tau^2$. We assume $g_1 \neq 0$ and $g_2 \neq 0$. We know from Bezout's lemma that if $a,~b\in \mathbb{Z}$ and $d=gcd(a,b)$, then there exists $\alpha_0$ and $\beta_0$ integers such that $d=\alpha_0a+\beta_0b$. All the other Bezout's coefficients of $d$ are given by 
$
\begin{cases}
\alpha_k= \alpha_0+k\frac{b}{d}\\
\beta_k= \beta_0-k\frac{a}{d}
\end{cases}
$
$k\in \mathbb{Z}$\\
Let's come back to our Drinfeld module $\phi_T=\gamma(T)+g_1\tau + g_2\tau^2$.\\ $d=gcd(q-1, q^2-1)=q-1$.\\
$d=q-1=-q(q-1)+1(q^2-1)$. The complete list of Bezout's coefficients of $d$ is given by:
$
\begin{cases}
\alpha_k= -q+k(q+1)=(k-1)q+k\\
\beta_k= 1-k
\end{cases}
$
$k \in \mathbb{Z}$.\\
Therefore the fine isomorphy invariant of $\phi$ is given by $$FI(\phi)=\left(g_1^{(k-1)q+k}.g_2^{1-k}~\left(modL^{*(q-1)} \right)\right)_{k \in \mathbb{Z}}$$
\end{ex}

\begin{df}{\label{JInv}}\cite[J-Invariants]{potemine1998minimal}~\\
Let $(k_1, \cdots, k_l)$ be a tuple with $1 \leq k_1< \cdots < k_l \leq r-1 $ and $\delta_1, \cdots, \delta_l$ be integers such that 
\begin{enumerate}[a)]
\item $\delta_1(q^{k_1}-1)+ \cdots +\delta_l(q^{k_l}-1)=\delta_r(q^r-1)$.
\item $0 \leq \delta_i \leq \frac{q^r-1}{q^{gcd(i,r)}-1}$~. for $i=1, \cdots, l.$
\item $gcd(\delta_1, \cdots, \delta_l, \delta_r)=1$
\end{enumerate}
The so-called basic J-invariants of the Drinfeld module $\phi$ are defined as $$J_{k_1 \cdots k_l}^{\delta_1 \cdots \delta_l}\left( \phi \right)=\frac{g_{k_1}^{\delta_1} \cdots g_{k_l}^{\delta_l}}{g_r^{\delta_r}}$$
\end{df}

\begin{thm}{\label{thmIso}} We keep the same notation above and we consider $\phi$ and $\psi: A \longrightarrow L\{\tau\}$ as two rank $r$ Drinfeld $A$-modules defined by $$ \phi_T= \gamma(T)+ g_1\tau + \cdots + g_r\tau^r \text{ and }  \psi_T= \gamma(T)+ g_1'\tau + \cdots + g_r'\tau^r$$.
The followings are equivalent 
\begin{enumerate}[(i)]
\item $\phi \overset{L}{\cong} \psi$
\item $\phi \overset{L^{sep}}{\cong} \psi \text{ and } \exists \lambda \in B,~~ FI_{\lambda}(\phi)=FI_{\lambda}(\psi)$
\item $\phi \overset{L^{sep}}{\cong} \psi \text{ and } FI(\phi)=FI(\psi)$
\end{enumerate}
\end{thm}
\begin{pf}
%We want to prove that $ \phi \overset{L}{\cong} \psi \Leftrightarrow \phi \overset{L^{sep}}{\cong} \psi \text{ and } FI(\phi)=FI(\psi) $\\
%\begin{itemize}
Our plan is to prove following the loop $(iii) \Rightarrow (ii) \Rightarrow (i) \Rightarrow (iii)$.\\
Let's assume $(iii)$. It obviously implies $(ii)$ since $B \neq \emptyset$.\\
Let's now assume for the second part of the proof that $\phi \overset{L^{sep}}{\cong} \psi$ and \\ $\exists \lambda=(\lambda_k)_{k \in I} \in B \text{ such that } FI_{\lambda}(\phi)=FI_{\lambda}(\psi)$.\\
We want to show that $\phi \overset{L}{\cong} \psi$.\\
$\phi \overset{L^{sep}}{\cong} \psi$ implies that there exists $x \in L^{sep} \text{ such that } \psi_T= x^{-1} \phi_T x$.\\
That is 
\begin{equation}\label{eq1}
\text{ for all } ~k\in I, ~ g_k'=g_k x^{q^k-1} 
\end{equation} 
$FI_{\lambda}(\phi)=FI_{\lambda}(\psi) \text{ implies } \displaystyle\prod_{k \in I}g_k'^{\lambda_k}=\prod_{k \in I}g_k^{\lambda_k}~mod L^{*d}$. That is
\begin{equation}\label{eq2}
\text{ there is} ~ y \in L^* \text{ such that } \displaystyle\prod_{k \in I}g_k'^{\lambda_k}=\prod_{k \in I}g_k^{\lambda_k} . y^d.
\end{equation}
From equation (\ref{eq1})  we get $g_k'^{\lambda_k}=g_k^{\lambda_k} x^{\lambda_k(q^k-1)}~~ \text{for all} ~k \in I  $\\ Thus 
\begin{equation}\label{eq3}
\displaystyle\prod_{k \in I}g_k'^{\lambda_k}=\prod_{k \in I}g_k^{\lambda_k} . x^{\sum\limits_{k \in I}\lambda_k(q^k-1)} =\prod_{k \in I}g_k^{\lambda_k} . x^d
\end{equation}
The equations (\ref{eq2}) and (\ref{eq3}) imply that $x^d=y^d$.\\
But $d=gcd(q^k-1, ~k\in I). \text{ That is for all } k\in I,\text{ there exists } \alpha_k \in \mathbb{Z} \text{ such that }$\\ $q^k-1=\alpha_k d$.\\ 
$\text{Hence } ~~x^{q^k-1}=x^{\alpha_k d}=\left(x^d\right)^{\alpha_k}=\left(y^d\right)^{\alpha_k}=y^{\alpha_k d}=y^{q^k-1}$.\\
Thus $\forall ~ k\in I~~g_k'=g_kx^{q^k-1}=g_ky^{q^k-1}$.\\
Therefore $\psi_T= y^{-1} \phi_T y$ and $y \in L^*$.\\
Hence $\phi \overset{L}{\cong} \psi$\\
For the last part of the proof we consider $(i)$. That is $\phi \overset{L}{\cong} \psi$. It obviously implies also that $\phi \overset{L^{sep}}{\cong} \psi$.\\Let's now check that $FI(\phi)=FI(\psi)$.\\ 
$\phi \overset{L}{\cong} \psi$ implies that there exists $x \in L$ such that $\psi_T=x^{-1} \phi_T x$.\\ That is, for all $k \in I, ~ g_k'=g_k x^{q^k-1}$. From The Bezout's lemma $B \neq \emptyset$. Let's then pick any $\lambda=(\lambda_k)_{k \in I} \in B$. We have $g_k'^{\lambda_k}=g_k^{\lambda_k}x^{\lambda_k(q^k-1)}$.\\
Thus $$\displaystyle\prod_{k \in I}g_k'^{\lambda_k}=\prod_{k \in I}g_k^{\lambda_k}\prod_{k \in I}x^{\lambda_k(q^k-1)}=\prod_{k \in I}g_k^{\lambda_k}.x^{\sum\limits_{k \in I}\lambda_k(q^k-1)}=\prod_{k \in I}g_k^{\lambda_k}. x^d$$
Therefore $\displaystyle\prod_{k \in I}g_k'^{\lambda_k}=\prod_{k \in I}g_k^{\lambda_k}. x^d, ~x \in L^*$. \\ Which implies $\displaystyle\prod_{k \in I}g_k'^{\lambda_k}=\prod_{k \in I}g_k^{\lambda_k}~mod L^{*d}$\\
Hence $FI_{\lambda}(\phi)= FI_{\lambda}(\psi)$.\\
Since $\lambda$ has been picked randomly, we can conclude that\\ $FI_{\lambda}(\phi)= FI_{\lambda}(\psi)~~\forall \lambda \in B$.\\
Therefore $FI(\phi)= FI(\psi)$. 
\end{pf}
\begin{rmk}
In the sequel, we might at some point abuse the language by considering as fine isomorphy invariants of $\phi$,  $FI_{\lambda_0}(\phi)\equiv FI(\phi)$ for some $\lambda_0 \in B$. As we can notice from the theorem above, this will not have any impact on the generality.
\end{rmk}

\begin{rmk} Potemine proved in \cite[Theorem 2.2]{potemine1998minimal} that
$$\phi \overset{L^{sep}}{\cong} \psi \Leftrightarrow J_{k_1 \cdots k_l}^{\delta_1 \cdots \delta_l}\left( \phi \right)=J_{k_1 \cdots k_l}^{\delta_1 \cdots \delta_l}\left( \psi \right) \text{ for any } (k_1, \cdots, k_l) \text{ and } (\delta_1, \cdots, \delta_l)$$
 $\text{ as defined above.}$\\
Taking it into account, one can reformulate the theorem \ref{thmIso} as follows.
\end{rmk}
\begin{thm}

$$\phi \overset{L}{\cong} \psi \Leftrightarrow J_{k_1 \cdots k_l}^{\delta_1 \cdots \delta_l}\left( \phi \right)=J_{k_1 \cdots k_l}^{\delta_1 \cdots \delta_l}\left( \psi \right) \text{ and } FI(\phi)=FI(\psi)$$
In other words, L-isomorphism classes of Drinfeld modules defined over the finite $A$-field L are determined by their fine isomorphy invariants and J-invariants.
\end{thm}
%\begin{pf}
%Straight forward from the theorem \ref{thmIso} above.
%\end{pf}
\begin{ex}
For the case of rank 2 Drinfeld $A$-modules, the only basic J-invariant is $J_1^{q+1}=\frac{g_1^{q+1}}{g_2}$. Here $d=$
$
\begin{cases}
\gcd(q-1, q^2-1)=q-1 & if ~g_1 \neq 0\\
q^2-1 & if~ g_1=0
\end{cases} 
$
Therefore $\lambda_1=$
$
\begin{cases}
-q & if~ g_1 \neq 0\\
0 & if~ g_1=0
\end{cases} 
$
and $\lambda_2=1$ in any case.\\
Thus $FI(\phi)=$
$
\begin{cases}
g_1^{-q}g_2\mod L^{*q-1} & if~ g_1 \neq 0\\
g_2 \mod L^{*q^2-1} & if~ g_1=0
\end{cases} 
$\\
The invariants $J_1^{q+1}$ and $FI(\phi)$ match clearly with the invariants describing the isomorphism classes of a rank 2 Dinfeld module as shown by Gekeler in \cite{gekeler2008frobenius}.

\end{ex}
\begin{ex}
Let's consider a rank 3 Drinfeld $A$-module defined over the field $L=\mathbb{F}_{25}=\mathbb{F}_5(\alpha)$ with $\alpha^2+4\alpha+2=0$. We take $A=\mathbb{F}_5[T]$. $L$ is an $A$-field defined by the ring homomorphism $\gamma :A \longrightarrow L, ~ T \mapsto \alpha $. \\
Let $\phi_T= \alpha + g_1\tau + g_2\tau^2+ g_3 \tau^3$.\\
Following the definition \ref{JInv}, one can easily compute the basic J-invariants of $\phi$ which are:\\
$J_{1,2}^{31,0}\left( \phi \right),~J_{1,2}^{1,5}\left( \phi \right),~J_{1,2}^{7,4}\left( \phi \right),~J_{1,2}^{8,9}\left( \phi \right),~J_{1,2}^{9,14}\left( \phi \right),~J_{1,2}^{10,19}\left( \phi \right),~J_{1,2}^{11,24}\left( \phi \right),~J_{1,2}^{12,29}\left( \phi \right)$ \\ $J_{1,2}^{13,3}\left( \phi \right),~J_{1,2}^{15,13}\left( \phi \right),~J_{1,2}^{17,23}\left( \phi \right),~J_{1,2}^{19,2}\left( \phi \right),~J_{1,2}^{20,7}\left( \phi \right),~J_{1,2}^{22,17}\left( \phi \right),~J_{1,2}^{23,22}\left( \phi \right),~J_{1,2}^{25,1}\left( \phi \right)$ \\$J_{1,2}^{27,11}\left( \phi \right),~J_{1,2}^{29,21}\left( \phi \right),~J_{1,2}^{31,31}\left( \phi \right),~$\\

\noindent The fine isomorphy invariant of $\phi$ is given by \\
$FI(\phi)=$
$ 
\begin{cases}
g_1 ~mod L^{*4} & \mbox{if } g_1 \neq 0 \\
\frac{g_3}{g_2^5}~mod L^{*4} & \mbox{if } g_1=0 \mbox{ and } g_2 \neq 0\\
g_3~mod L^{*124} & \mbox{if } g_1=g_2=0
\end{cases}
$\\

\noindent Therefore the isomorphism class of $\phi$ is parametrized by those 20 invariants
\end{ex}
%\begin{rmk}
%For Drinfeld modules of a given rank defined over a finite field, Potemine proved in \cite{potemine1998minimal} that the number of isomorphism classes is given by:
%$$ \#Cl(D^r/L)= q^{n_r}-1 + \displaystyle\sum_{(i_1, \cdots, i_s) \in I_0}\left(q^{gcd(i_1, \cdots,i_s, n_r)}-1 \right)\left(q^n-1\right)^s+\left(q-1\right)\left[q^{(r-1)n}-q^{r-\varphi(n_r, r)n}\right]$$
%Where $n=[L: \mathbb{F}_q],~ n_r=gcd(n, r).~~I_0$ is the power set of $\{1, \cdots, r\}$ and $I_1 \subset I_0$ is made up of subsets $(i_1, \cdots , i_s)$ such that $gcd(i_1, n_r)> 1, \cdots, gcd(i_s, n_r)> 1.~~\varphi(n_r, r)$ is the number of integers $< r$ and coprime with $n_r$. 
%\end{rmk}
\begin{rmk}
Each isomorphism class has a finite number of elements.\\ Indeed $\#Cl(\phi) \leq \# L^{*d}$.
\end{rmk}
\noindent We provide in the sequel an algorithm generating the isomorphism classes of rank $r$ Drinfeld modules in a given isogeny class.
\begin{algo}\label{algo_iso}[Isomorphism classes of a Drinfeld modules]~\\
\textbf{Inputs}: Isogeny class defined by $M(x)= x^{r_1} + a_1x^{r_1-1} + \cdots + a_{r_1-1}x + \mu \mathfrak{p}_v^{\frac{m}{r_2}}$.\\
\textbf{Ouputs}: Isomorphism classes of Drinfeld modules in the isogeny class defined by $M(x)$ 

\begin{enumerate}[1-]
\item Set $\phi_T= g_r\tau^r + \cdots + g_1\tau + \gamma(T)$ and solve the equation (system of equations) given by $\tau^{sr_1}+ a_1(\phi_T)\tau^{s(r_1-1)} + \cdots + a_{r_1-1}(\phi_T)\tau^s + \mu \mathfrak{p}_v(\phi_T)^{\frac{m}{r_2}}=0$.\\ 
Where $s=[L: \mathbb{F}_q]$. Let $\Gamma$ be the set of all solutions of that equation. 
\item Pick a Drinfeld module $\phi \in \Gamma$. We assume $\phi_T= g_r\tau^r + \cdots + g_1\tau + \gamma(T)$.
\item Compute the fine isomorphy invariant and the $J$-invariants of $\phi$. i.e. $FI(\phi)$ and $J_{k_1 \cdots k_l}^{\delta_1 \cdots \delta_l}\left(\phi \right).$
\item for $\psi$ in $\Gamma$:
	Compute $FI(\psi)$ and $J_{k_1 \cdots k_l}^{\delta_1 \cdots \delta_l}\left(\psi \right).$\\
	If $FI(\psi) = FI(\phi)$ and  $J_{k_1 \cdots k_l}^{\delta_1 \cdots \delta_l}\left(\psi \right) =J_{k_1 \cdots k_l}^{\delta_1 \cdots \delta_l}\left(\phi \right)$:\\
	Then store $\psi$ in the isomorphism class of $\phi$.\\
\item Pick another $\phi$ in $\Gamma$ which is not in the previously computed isomorphism classes and move to step 3.
\item If the set $\Gamma$ is exhausted then output the isomorphism classes and exit.

%\item Compute $d=gcd(q^k-1,~k\in I)$ and the Bezout coefficients $\lambda_k$ such that $d=\displaystyle\sum_{k \in I} \lambda_k(q^k-1)$ where $I=\{i=1, \cdots, r;~g_i \neq 0\}$
%\item Choose $y=x_0^d \in L^{*d}$
%\item Compute $g_k'=g_k . x_0^{q^k-1}~~ \forall k \in I$
%\item Compute the basic J-invariants $J_{k_1 \cdots k_l}^{\delta_1 \cdots \delta_l}\left(\phi \right)$ of $\phi$.
%\item Consider the Drinfeld module $\psi$ defined by $\psi_T=\gamma(T)+ \displaystyle\sum_{k \in I}g_k'\tau^k$
%\item Compute the basic J-invariants $J_{k_1 \cdots k_l}^{\delta_1 \cdots \delta_l}\left(\psi \right)$ of $\psi$.
%\item If $J_{k_1 \cdots k_l}^{\delta_1 \cdots \delta_l}\left(\psi \right)=J_{k_1 \cdots k_l}^{\delta_1 \cdots \delta_l}\left(\phi \right)$ then store $\psi$ in the isomorphism class of $\phi$ and go back to step 3.\\Else go back to step 3.
\end{enumerate} 
%The algorithm stops when we have exhausted all the elements of $L^{*d}$ 

\end{algo}
%\begin{rmk}
%This algorithm also works for any isogeny class defined by a Weil polynomial of the form 
%$$M(x)= x^{r_1} + a_1x^{r_1-1} + \cdots + a_{r_1-1}x + \mu Q^{1/r_2} \text{ with } r=r_1r_2 \text{ and } r_2 \mid m.$$
%\end{rmk}

\section{Application: more specific description for the case of rank 3 Drinfeld modules}

\subsection{Isogeny classes for rank 3 Drinfeld modules}
We keep the same data as before. That is $A=\mathbb{F}_q[T], ~ k=\mathbb{F}_q(T)$ with a distinguished place at infinity $\infty$.\\
As we have seen before, the isogeny classes are given by the following rank 3 Weil polynomials:
\begin{itemize}
\item $M(x)=x^3+a_1x^2 + a_2x + \mu \mathfrak{p}_v^m \in A[x]$ with $\mu \in \mathbb{F}_q$. Where $\deg a_1 \leq \frac{m\deg\mathfrak{p}_v }{3}$ and $\deg a_2 \leq \frac{2m\deg\mathfrak{p}_v}{3}$ such that the resultant modulo $\mathfrak{p}_v$ of any two irreducible factors of $M(x) \mod \mathfrak{p}_v^n$ is non-zero and\\ $\overline{M_0(x)} \equiv x^3 + \frac{a_1}{T^s}x^2 + \frac{a_2}{T^{2s}}x + \mu\frac{\mathfrak{p}_v^m}{T^{3s}} \mod \frac{1}{T^h}$ is irreducible. \\
Where $h=v_{\infty}\left(disc\left(M(x)\right)\right) + sr(r-1) +1$ and $n=v\left(disc\left(M(x)\right)\right) +1$ (see algorithm~\ref{algo}).
\item $M(x)= x-\mu \mathfrak{p}_v^{\frac{m}{3}}$ with $3 \vert m$ and $\mu \in \mathbb{F}_q^*$
\end{itemize} 
We provide in the sequel some more specific results that help to quickly identify rank 3 Weil polynomials by more or less just looking atthe ``size" of the coefficient of the polynomials.\\
Before that let us recall the notion of standard form of a cubic polynomial.

\begin{df}[Standard form] \label{def_stdform}~\\
Let $k(\tilde{\pi})/k$ be a cubic function field. The minimal polynomial $M_0(x) \in A[x]$ of $\tilde{\pi}$ is said to be in the standard form if $M_0(x)=x^3+ax+b$ with $a$ and $b \in A$ satisfying the following: $$\text{There is no } c \in A \text{ such that } c^2 \vert a \text{ and } c^3 \vert b.$$
\end{df}
\begin{rmk}
Let $M(x)=x^3+a_1x^2+a_2x +\mu \mathfrak{p}_v^m$ be a potential Weil polynomial whose corresponding cubic field is $k(\pi)/k$.\\
If $char(k) \neq 3$, setting $x=y-\frac{a_1}{3}$, one can transform $$M(x)=x^3 + a_1x^2 + a_2x + \mu \mathfrak{p}_v^m$$ into a polynomial of the form $$y^3+b_1y + b_2 \in A[y] \text{ where } b_1=\frac{-a_1^2}{3}+a_2, ~ b_2=\frac{2a_1^3}{27} -\frac{a_1a_2}{3} + \mu Q.$$
 One can therefore convert the polynomial $N(y)=y^3+b_1y +b_2 \in A[y]$ into a standard polynomial $x^3 +c_1x +c_2$. By ``converting" we mean getting from the irreducible polynomial $y^3+b_1y + b_2$ an irreducible polynomial in the standard form $M_0(x)=x^3 +c_1x +c_2$ whose any root $\tilde{\pi}$ is such that $k(\tilde{\pi}) \backsimeq k(\pi)$ (i.e. $k(\tilde{\pi})$ and $k(\pi)$ define the same function field).\\
In fact doing it, is really a simple exercise. One takes the square-free factorizations of $b_1$ and $b_2$. That is $b_1 = \mu_1\displaystyle\prod_{i=1}^{n_1} b_{1i}^i$ and $b_2 = \mu_2\displaystyle\prod_{j=1}^{n_2} b_{2j}^j$ where $\mu_1, ~\mu_2 \in \mathbb{F}_q$ and $b_{1i}~ i=1, \cdots, n_1$ (resp. $b_{2j}~ j=1, \cdots, n_2$) are pairwise coprime square-free elements of $A$. We set $g_1 = \displaystyle\prod_{i=1}^{n_1} b_{1i}^{\floor{\frac{i}{2}}}$ and $g_2 = \displaystyle\prod_{j=1}^{n_2} b_{2j}^{\floor{\frac{j}{3}}}$. Taking $c_1=\frac{b_1}{\gcd(g_1, g_2)^2}$ and $c_2=\frac{b_2}{\gcd(g_1, g_2)^3}$,  we have that $M_0(x)=x^3 +c_1x +c_2$ is a polynomial in the standard form in $A[x]$. In addition we have the following:\\ 
$\pi$ is a root of $M(x)$ if and only if $\pi+\frac{a_1}{3}$ is a root of $N(y)=y^3+b_1y +b_2$ if and only if $\tilde{\pi}=\frac{\pi + \frac{a_1}{3}}{\gcd(g_1, g_2)}$ is a root of $M_0(x)= x^3 +c_1x +c_2$.\\
Therefore $disc\left( M(x) \right) = disc\left(N(y)\right)$ and $ind(\pi)= ind\left(\pi + \frac{a_1}{3}\right)$.\\ But $ind(\tilde{\pi})=\frac{ind(\pi)}{\gcd(g_1, g_2)^3}$ because $disc\left( M_0(x)\right)= \frac{disc\left(M(x)\right)}{\gcd(g_1, g_2)^6}$.\\ Also, $k(\pi)= k\left(\pi + \frac{a_1}{3}\right)=k(\tilde{\pi}).$
\end{rmk}
%As a consequence of proposition~\ref{prop2} in the special case of the degree $3$ polynomial $M_0(x)$ in the standard form, we have the following:
%one gets that $M_0(x)=x^3 +c_1x +c_2 \in A[x]$. In addition $\pi$ is a root of $Y^3 + b_1Y + b_2$ if and only if $\tilde{\pi}=\frac{\pi}{\gcd(g_1, g_2)}$ is a root of $M_0(x)$. That means $k(\tilde{\pi}) \backsimeq k(\pi).$ Also $disc\left(M_0(x) \right) = \frac{disc\left(M(x)\right)}{\gcd(g_1, g_2)^6}$.\\   
%The algorithm~\ref{algo} can be speeded up using the following propositions:\\
%Let $k(\pi)/k$ be our cubic function field such that the standard form of the minimal polynomial of $\pi$ is given by $$M_0(x)=x^3+c_1x +c_2$$
\begin{prop} \cite[theorem 4.2]{scheidler2004algorithmic}\label{propos1}~\\
Let $M_0(x)=x^3 +c_1x +c_2$ be the standard form of the minimal polynomial $M(x)$ of $\pi$.\\
There is a unique place of $k(\pi)$ above the place at infinity $\infty$ of $k$ only in the following cases.
\begin{itemize}
\item[$(s1)$] $3\deg c_1 < 2\deg c_2, ~ \deg c_2 \equiv 0 \mod 3$ and $LC(c_2)$ is not a cube in $\mathbb{F}_q$. $LC(?)$ denotes here the leading coefficient of the argument.
%\item[$(s2)$] $3\deg c_1 = 2\deg c_2, ~ \deg\left((disc\left(M_0(x)\right)\right)$ is even, $LC\left(disc\left(M_0(x)\right)\right)$ is a square in $\mathbb{F}_q$ and $\delta=-\frac{3}{2}LC\left(9c_2 + \sqrt{-3disc\left(M_0(x)\right)}\right)$ is not a cube in $\mathbb{F}_q(\xi)$. Where $\xi$ denotes a primitive cubic root of unity.
\item[$(s2)$] $3\deg c_1 = 2\deg c_2,~ 4LC(c_1)^3+27LC(c_2)^2 \neq 0$ and\\ $x^3+ LC(c_1)x + LC(c_2)$ has no root in $\mathbb{F}_q$.
\item[$(s3)$] $3\deg c_1 < 2\deg c_2 \text{ and } \deg c_2 \not\equiv 0 \mod 3 $
\end{itemize}
%One can notice that $\deg\left(disc\left(M_0(x)\right)\right)$ is even in each of those cases.
\end{prop}

\begin{prop}\label{c2-rank3}
Let $M(x)=x^3+a_1x^2+a_2x + \mu\mathfrak{p}_v^m \in A[x]$ be as mentioned before.
\begin{enumerate}
\item If $\mathfrak{p}_v \mid a_2$ and $\mathfrak{p}_v \nmid a_1$ then there is a unique zero of $\pi$ in $k(\pi)$ above the place $v$ if and only if $v(a_2) \geq \frac{m}{2}$.
\item If $\mathfrak{p}_v \mid a_2$ and $\mathfrak{p}_v \mid a_1$ then there is a unique zero of $\pi$ in $k(\pi)$ above the place $v$ if and only if there is a unique place of $k(\pi)$ above $v$ (i.e. if and only if $M(x)$ is irreducible over the completion field $k_v$).
\item If $\mathfrak{p}_v \nmid a_2$ then there is a unique  zero of $\pi$ in $k(\pi)$ above $v$. 
\end{enumerate}
\end{prop}

\noindent Before proving this proposition, let us recall the following lemma, known as Hensel lemma or Hensel lifting.
\begin{lem}\label{Hensel_lem}
Let $M(x) \in A[x]$ and $\mathfrak{p}$ be a prime in $A$. Let $m, n \in \mathbb{N}$ with $m \leq n$
\begin{itemize}
\item If $M(x_0) \equiv 0 \mod \mathfrak{p}^n$ and $M'(x_0) \not\equiv 0 \mod \mathfrak{p}$ then there exists a unique lifting of $x_0$ modulo $\mathfrak{p}^{n+m}$. i.e. there exists a unique $x_1 \in A$ such that $M(x_1) \equiv 0 \mod \mathfrak{p}^{n+m}$ and $x_1 \equiv x_0 \mod \mathfrak{p}^n$.
\item If  $M(x_0) \equiv 0 \mod \mathfrak{p}^n$ and $M'(x_0) \equiv 0 \mod \mathfrak{p}$ then we have two possibilities:
\begin{itemize}
\item If $M(x_0) \not\equiv 0 \mod \mathfrak{p}^{n+1}$ then there is no lifting of $x_0$ modulo $\mathfrak{p}^{n+1}$.
\item If $M(x_0) \equiv 0 \mod \mathfrak{p}^{n+1}$ then every lifting of $x_0$ modulo $\mathfrak{p}^{n+1}$ is a zero of $M(x)$ modulo $\mathfrak{p}^{n+1}$. 
\end{itemize}  
\end{itemize}
\end{lem}

\begin{pf}[Proof of proposition~\ref{c2-rank3}]~\\
\begin{enumerate}
\item We assume here that $\mathfrak{p}_v \mid a_2$ and $\mathfrak{p}_v \nmid a_1$.

\fbox{\parbox[b]{0.5cm}{$\Rightarrow$}} We assume that there is a unique zero of $\pi$ in $k(\pi)$ above $v$.

$M(x) \equiv x^2(x+a_1) \mod \mathfrak{p}_v$ and $\mathfrak{p}_v \nmid a_1$. That means $0$ (as double root) and $-a_1$ are the roots of $M(x)$ module $\mathfrak{p}_v$.\\ 
Using the Hensel lemma~\ref{Hensel_lem}, one can lift these roots modulo $\mathfrak{p}_v^l ~(\text{for } l\geq 1)$ as long as $M(0) \equiv 0 \mod \mathfrak{p}_v^l$. \\
We know that $disc\left( M(x)\right) = (a_1^2-4a_2)a_2^2 + \mathfrak{p}_v^m(-4a_1^3-27\mathfrak{p}_v^m + 18a_1a_2)$\\
Let us assume that $v(a_2) < \frac{m}{2}$.\\
That means $v(a_2^2) < m$. Since $\mathfrak{p}_v \nmid a_1$ and $\mathfrak{p}_v \mid a_2$, $v(a_1^2-4a_2)=0$ and $v(-4a_1^3-27\mathfrak{p}_v^m + 18a_1a_2)=0$. In other word 
$$v\left( disc\left(M(x)\right)\right) = v(a_2^2) < m.$$
For any $n \in \mathbb{N}$ with $n \leq m$, $M(0) \equiv 0 \mod \mathfrak{p}_v^n$. One can therefore lift the root $x_0=0$ modulo $\mathfrak{p}_v$ to roots modulo $\mathfrak{p}_v^n$ for $n=v(a_2^2) +1$ and the (simple) root $x_1=-a_1$ modulo $\mathfrak{p}_v$  to a root modulo $\mathfrak{p}_v^n$. One gets then
$$M(x) \equiv M_1(x) \cdot M_2(x) \cdot M_3(x) \mod \mathfrak{p}_v^{v\left(disc\left(M(x)\right)\right)+1} $$
With $M_1(x) \equiv M_2(x) \equiv x \mod \mathfrak{p}_v$ and $M_3(x) \equiv x+a_1 \mod \mathfrak{p}_v$.\\
Thus $Res\left(M_1(x), M_2(x)\right) \equiv 0 \mod \mathfrak{p}_v$ which contradicts the fact that there is a unique zero of $\pi$ in $k(\pi)$ above $v$ (see lemma~\ref{prop6} and corollary ~\ref{cor1}).\\
Therefore $v(a_2) \geq \frac{m}{2}.$

\fbox{\parbox[b]{0.5cm}{$\Leftarrow$}} Let us assume conversely that $v(a_2) \geq \frac{m}{2}$. We want to show that there is a unique zero of $\pi$ in $k(\pi)$ above $v$.\\
We recall that $disc\left( M(x)\right) = (a_1^2-4a_2)a_2^2 + \mathfrak{p}_v^m(-4a_1^3-27\mathfrak{p}_v^m + 18a_1a_2)$.\\
$\mathfrak{p}_v \mid a_2$ and $\mathfrak{p}_v \nmid a_1$ implies that $v(a_1^2-4a_2)=v(-4a_1^3-27\mathfrak{p}_v^m + 18a_1a_2)=0.$ In addition, $v(a_2^2)=2v(a_2) \geq m$. Thus $v\left(disc\left(M(x)\right)\right) \geq m.$\\
But $M(x) \equiv x^2(x+a_1) \mod \mathfrak{p}_v$ with $\mathfrak{p}_v \nmid a_1$.\\
The root $x_0=0$ of $M(x) \mod \mathfrak{p}_v$ can be lifted to a root of $M(x) \mod \mathfrak{p}_v^n$ for $n\leq m$. But since for $n \geq m+1$ $M(0) \not\equiv 0 \mod \mathfrak{p}_v^n$, there is no lifting of $x_0$ to a root of $M(x) \mod \mathfrak{p}_v^n$ (see Hensel lemma~\ref{Hensel_lem}). In other words, we cannot have $M(x) \equiv M_1(x) \cdot M_2(x) \cdot M_3(x) \mod \mathfrak{p}_v^{v\left(disc\left(M(x)\right)\right)+1}$ with $M_1(x) \equiv M_2(x) \equiv x \mod \mathfrak{p}_v$ and $M_3(x) \equiv x+a_1 \mod \mathfrak{p}_v$.\\
Therefore we are only left with the possibility\\ 
$M(x) \equiv M_1(x) \cdot M_2(x) \mod \mathfrak{p}_v^{v\left(disc\left(M(x)\right)\right)+1} $ with\\
$M_1(x) \equiv x^2 \mod \mathfrak{p}_v$ and $M_2(x) \equiv x+a_1 \mod \mathfrak{p}_v$ (see~\cite[Corollary 2.4]{von1996factorization}). We therefore clearly have $Res\left(M_1(x), M_2(x)\right) \not\equiv 0 \mod \mathfrak{p}_v$ since $\mathfrak{p}_v \nmid a_1$.\\
Hence there is a unique zero of $\pi$ in $k(\pi)$ above the place $v$.
\item we assume here that $\mathfrak{p}_v \mid a_1$ and $\mathfrak{p}_v \mid a_2$.\\
$M(\pi)=0$ implies that $\pi^3=-a_1\pi^2-a_2\pi-\mu\mathfrak{p}_v^m= \mathfrak{p}_v\left( -b_1\pi^2-b_2\pi - \mu\mathfrak{p}_v^{m-1} \right)$ where $a_i=b_i\cdot \mathfrak{p}_v$. In other words $\mathfrak{p}_v$ divides $\pi$. That means any place of $k(\pi)$ above $v$ is a zero of $\pi$.\\
Therefore there is a unique zero of $\pi$ in $k(\pi)$ above $v$ if and only if there is a unique place of $k(\pi)$ above $v$.  
\item This case has already been shown in proposition~\ref{c2-cor}.
\end{enumerate}
\end{pf}

\noindent We summarize our previous results in the following theorem.

\begin{thm}
Let $M(x)=x^3 + a_1x^2 + a_2x + \mu\mathfrak{p}_v^m \in A[x]$ be a potential Weil polynomial. i.e. $\deg a_i \leq \frac{im\deg \mathfrak{p}_v}{3}$ and $M(x)$ irreducible over $k$. We also consider $M_0(x)=x^3 + c_1x + c_2$ the standard form of $M(x)$.
\begin{enumerate}
\item There is a unique place of $k(\pi)$ lying over the place at $\infty$ of $k$ if and only if one of the following holds.
\begin{itemize}
\item[$(s1)$] $3\deg c_1 < 2\deg c_2, ~ \deg c_2 \equiv 0 \mod 3$ and $LC(c_2)$ is not a cube in $\mathbb{F}_q$. 
\item[$(s2)$] $3\deg c_1 = 2\deg c_2,~ 4LC(c_1)^3+27LC(c_2)^2 \neq 0$ and\\ $x^3+ LC(c_1)x + LC(c_2)$ has no root in $\mathbb{F}_q$.
\item[$(s3)$] $3\deg c_1 < 2\deg c_2 \text{ and } \deg c_2 \not\equiv 0 \mod 3 $\\
$LC(?)$ denotes here the leading coefficient of the argument.
\end{itemize}
\item There is a unique zero of $\pi$ in $k(\pi)$ lying over the place $v$ of $k$ if and only if one of the following holds.
\begin{itemize}
\item[$(s4)$] $\mathfrak{p}_v \mid a_2, ~~ \mathfrak{p}_v \nmid a_1$ and $v(a_2) \geq \frac{m}{2}$
\item[$(s5)$] $\mathfrak{p}_v \mid a_2, ~~ \mathfrak{p}_v \mid a_1$ and $M(x) \mod \mathfrak{p}_v^n$ is irreducible.\\
Where $n=v\left(disc\left(M(x)\right)\right) + 1$.
\item[$(s6)$] $\mathfrak{p}_v \nmid a_2$.
\end{itemize}
\end{enumerate}
\end{thm}

\noindent Using the previous results, one can therefore get a more specific version of the algorithm \ref{algo} for $r=3$ as follows:
\begin{algo}\label{newalgo}
\textbf{Input}: $M(x)=x^3+a_1x^2 +a_2x +\mu Q \in A[x]$ irreducible polynomial defining the cubic field $k(\pi)/k$.\\
\textbf{Ouput}: \textbf{True} if $M(x)$ is a Weil polynomial and \textbf{False} otherwise.
\begin{enumerate}
\item Compute $b_1=\frac{-a_1^2}{3}+a_2; ~~b_2=\frac{2a_1^3}{27} - \frac{a_1a_2}{3} + \mu Q$.
\item Compute the square-free decomposition of $b_1$ and $b_2$:\\
$b_1=\mu_1\displaystyle\prod_{i=1}^{n_1}b_{1i}^i,~~
b_2=\mu_2\displaystyle\prod_{j=1}^{n_2}b_{1j}^j$\\
Set $g_1=\displaystyle\prod_{i=1}^{n_1}b_{1i}^{\floor{\frac{i}{2}}} \text{ and }
g_2=\displaystyle\prod_{j=1}^{n_2}b_{1j}^{\floor{\frac{j}{3}}}$
\item Compute $c_1=\frac{b_1}{gcd(g_1, g_2)^2}$ and $c_2=\frac{b_2}{gcd(g_1, g_2)^3}$
\item If $c_1$ and $c_2$ fulfill one of the statements $(s1)$, $(s2)$ or $(s3)$ of proposition~\ref{propos1} then move to the next step. Otherwise output \textbf{False} and exit
\item Compute $n=v\left(disc(M(x))\right) + 1$ and\\ $\overline{M(x)} \equiv x^3+ a_1x^2 + a_2x + \mu Q \mod \mathfrak{p}_v^n$.\\ 
If $\mathfrak{p}_v \mid a_2$ and $\mathfrak{p}_v \nmid a_1$ and $v(a_2) \geq \frac{m}{2}$ then output \textbf{True} and exit.\\
Else if $\mathfrak{p}_v \mid a_2$ and $\mathfrak{p}_v \mid a_1$ and $\overline{M(x)}$ is irreducible then the output \textbf{True} and exit.\\
Else if $\mathfrak{p}_v \nmid a_2$ then output \textbf{True} and exit.\\
Else output \textbf{False} and exit. 
\end{enumerate}
\end{algo}

\subsection{Example of computation of isomorphism classes in a rank 3 isogeny class}

Here we mainly explain how the computation can be done and we provide a concrete example.\\
We consider the isogeny class defined by the polynomial $$M(x)=x^3+a_1(T)x^2 + a_2(T)x + \mu Q(T)$$
We want to list all the isomorphism classes of Drinfeld modules in this isogeny class. We know that the Frobenius endomorphism $\pi= \tau^s$ (with $s=[L: \mathbb{F}_q]$) is a root of $M(x)$. That means $$\tau^{3s} + a_1(T)\tau^{2s} + a_2(T)\tau^{s} + \mu Q(T)=0.$$
By definition of the action of the Drinfeld module $\phi$ we have
$$\tau^{3s} + a_1(\phi_T)\tau^{2s} + a_2(\phi_T)\tau^{s} + \mu Q(\phi_T)=0 \hspace{2cm} (\star)$$
We consider $(\star)$ as an equation with unknown $\phi_T$. This equation can be solved by setting $\phi_T=\gamma(T) + \alpha_1\tau + \alpha_2\tau^2+ \alpha_3\tau^3.$ We recall that $\gamma(T)$ is already known since $\gamma$ is the ring homomorphism defining the A-field L. One can therefore plug $\phi_T$ in the equation $(\star)$ and get a non-linear system of equation (with unknowns $\alpha_{i's}$). Even though the system is non-linear, a way to solve it can be by "brute force". That is, looking for all tuples $(\alpha_1, \alpha_2, \alpha_3) \in L^3$ solutions of the system. Since $L$ is finite, we have finitely many such tuples. Each of those solutions yields a Drinfeld module $\phi$ defined by $\phi_T=\gamma(T) + \alpha_1\tau + \alpha_2\tau^2+ \alpha_3\tau^3.$ We therefore gather those Drinfeld modules with respect to their isomorphism classes by computing and comparing their $J$-invariants and fine isomorphy invariants.\\
Let us have a look at a concrete example.\\
Let $A=\mathbb{F}_5[T],~ k=\mathbb{F}_5(T), ~L=\mathbb{F}_5(\alpha)$ with $\alpha^2+4\alpha+2=0$. $L$ is an $A$-field defined by $\gamma : A \longrightarrow L, ~f(T) \longmapsto f(0)$. The $A$-characteristic of $L$ is $T$ because $\mathfrak{p}_v=Ker\gamma = T\cdot A$ is the ideal generated by $T$.\\ $m=[L: A/\mathfrak{p}_v]=[L: A/T\cdot A]=[\mathbb{F}_5(\alpha): \mathbb{F}_5]=2$. We consider the polynomial $$M(x)=x^3+ 3x^2 + (1+T)x + T^2$$
\underline{\textbf{Claim}}: $M(x)$ is a Weil polynomial.\\
first of all $M(x)$ is irreducible in $A[x]$ and therefore (Gauss lemma) is also irreducible in $k[x]$. One easily shows using the algorithm~\ref{algo} that
\begin{itemize}
\item $\overline{M_0(x)}=x^3 + \frac{3}{T}x^2 + \frac{1+T}{T^2}x + \frac{T^2}{T^3} \mod \frac{1}{T^3} \equiv x^3 + \frac{3}{T}x^2 + \frac{1+T}{T^2}x + \frac{1}{T} \mod \frac{1}{T^3}$ ($h=3$) is irreducible.
\item $\overline{M(x)}= x^3+ 3x^2 + (1+T)x + T^2 \mod T^2 \equiv x(x^2 + 3x + 1 + T) \mod T^2$ ($n=2$) and we clearly have $Res(x, x^2 + 3x + 1 + T ) \mod T \not\equiv 0$.
\end{itemize}
Hence $M(x)$ defines an isogeny class of Drinfeld modules.\\
We aim to list (as explained before) all the isomorphism classes of Drinfeld modules in the isogeny class defined by $M(x)=x^3+ 3x^2 + (1+T)x + T^2.$\\
$\pi = \tau^s$ with $s=[L: \mathbb{F}_5]=2$. i.e. $\pi= \tau^2$. In addition $M(\pi)=0$ i.e. 
$$\tau^6 + 3\tau^4 + (1+ \phi_T)\tau^2 + \phi_T^2=0$$
That means $\phi_T^2 + \phi_T \tau^2 + \tau^6 + 3\tau^4 + \tau^2 = 0$. We clearly see from the Weil polynomial that $T \in ker \gamma.$ i.e. $\gamma(T)=0.$\\
We can therefore set $\phi_T=\alpha_1\tau + \alpha_2 \tau^2 + \alpha_3\tau^3 \in L\{ \tau \}$.\\
i.e. $(\alpha_1\tau + \alpha_2 \tau^2 + \alpha_3\tau^3)^2 + ( \alpha_1\tau + \alpha_2 \tau^2 + \alpha_3\tau^3)\tau^2 + \tau^6 + 3\tau^4 + \tau^2 = 0$.\\
Solving this equation yields the following Drinfeld modules:\\
%$$\phi_T= 4\tau + 2\tau^2 + \tau^3 \text{ and } \psi_T = 4\tau + 2\tau^2 + 4\tau^3.$$
%%%%%%%%%%%%%%%%%%%%%%%%%%%%%%%%%%%%%%%%%%%%%%%%%%%%
%%%%%%%%%%%%%%%%%%%%%%%%%%%%%%%%%%%%%%%%%%%%%%%%%%%%
$$ 
\begin{array}{|cc|}
\hline
\hspace{7cm} \phi(T) & \\
\hline

(\alpha + 3)\tau + 2\tau^2 + (4\alpha + 4)\tau^3 &
(\alpha + 3)\tau + 2\tau^2 + 3\tau^3 \\
(\alpha + 3)\tau + (2\alpha + 1)\tau^2 + (\alpha + 3)\tau^3 & 
(\alpha + 3)\tau + 4\alpha\tau^2 + 2\tau^3 \\
(\alpha + 3)\tau + 4\alpha\tau^2 + (\alpha + 1)\tau^3 & 
(\alpha + 3)\tau + (3\alpha + 3)\tau^2 + (\alpha + 3)\tau^3 \\ 
(\alpha + 3)\tau + (\alpha + 4)\tau^2 + 2\tau^3 &
(\alpha + 3)\tau + (\alpha + 4)\tau^2 + (\alpha + 1)\tau^3 \\ 
2\tau + 2\tau^2 + (4\alpha + 2)\tau^3 &
2\tau + 2\tau^2 + (\alpha + 1)\tau^3 \\
2\tau + (2\alpha + 1)\tau^2 + 2\tau^3 & 
2\tau + 4\alpha\tau^2 + (\alpha + 3)\tau^3 \\ 
2\tau + 4\alpha\tau^2 + (4\alpha + 4)\tau^3 & 
2\tau + (3\alpha + 3)\tau^2 + 2\tau^3 \\
2\tau + (\alpha + 4)\tau^2 + (\alpha + 3)\tau^3 & 
2\tau + (\alpha + 4)\tau^2 + (4\alpha + 4)\tau^3 \\ 
(4\alpha + 4)\tau + 2\tau^2 + (\alpha + 3)\tau^3 &
(4\alpha + 4)\tau + 2\tau^2 + 3\tau^3 \\
(4\alpha + 4)\tau + (2\alpha + 1)\tau^2 + (4\alpha + 4)\tau^3 & 
(4\alpha + 4)\tau + 4\alpha\tau^2 + 2\tau^3 \\
(4\alpha + 4)\tau + 4\alpha\tau^2 + (4\alpha + 2)\tau^3 & 
(4\alpha + 4)\tau + (3\alpha + 3)\tau^2 + (4\alpha + 4)\tau^3 \\ 
(4\alpha + 4)\tau + (\alpha + 4)\tau^2 + 2\tau^3 &
(4\alpha + 4)\tau + (\alpha + 4)\tau^2 + (4\alpha + 2)\tau^3 \\ 
(4\alpha + 2)\tau + 2\tau^2 + 2\tau^3 &
(4\alpha + 2)\tau + 2\tau^2 + (\alpha + 1)\tau^3 \\ 
(4\alpha + 2)\tau + (2\alpha + 1)\tau^2 + (4\alpha + 2)\tau^3 & 
(4\alpha + 2)\tau + 4\alpha\tau^2 + (4\alpha + 4)\tau^3 \\
(4\alpha + 2)\tau + 4\alpha\tau^2 + 3\tau^3 &
(4\alpha + 2)\tau + (3\alpha + 3)\tau^2 + (4\alpha + 2)\tau^3 \\ 
(4\alpha + 2)\tau + (\alpha + 4)\tau^2 + (4\alpha + 4)\tau^3 &
(4\alpha + 2)\tau + (\alpha + 4)\tau^2 + 3\tau^3 \\
3\tau + 2\tau^2 + (\alpha + 3)\tau^3 &
3\tau + 2\tau^2 + (4\alpha + 4)\tau^3 \\ 
3\tau + (2\alpha + 1)\tau^2 + 3\tau^3 &
3\tau + 4\alpha\tau^2 + (4\alpha + 2)\tau^3 \\ 
3\tau + 4\alpha\tau^2 + (\alpha + 1)\tau^3 &
3\tau + (3\alpha + 3)\tau^2 + 3\tau^3 \\
3\tau + (\alpha + 4)\tau^2 + (4\alpha + 2)\tau^3 & 
3\tau + (\alpha + 4)\tau^2 + (\alpha + 1)\tau^3 \\
(\alpha + 1)\tau + 2\tau^2 + 2\tau^3 &
(\alpha + 1)\tau + 2\tau^2 + (4\alpha + 2)\tau^3 \\ 
(\alpha + 1)\tau + (2\alpha + 1)\tau^2 + (\alpha + 1)\tau^3 & 
(\alpha + 1)\tau + 4\alpha\tau^2 + (\alpha + 3)\tau^3 \\
(\alpha + 1)\tau + 4\alpha\tau^2 + 3\tau^3 &
(\alpha + 1)\tau + (3\alpha + 3)\tau^2 + (\alpha + 1)\tau^3 \\ 
(\alpha + 1)\tau + (\alpha + 4)\tau^2 + (\alpha + 3)\tau^3 &
(\alpha + 1)\tau + (\alpha + 4)\tau^2 + 3\tau^3 \\
\hline

\end{array} 
$$
 ~\\
We have implemented a SAGE code adapted to algorithm~\ref{algo_iso} in order to gather these Drinfeld modules with respect to their isomorphism classes and we got the following:\\
~\\
$
\begin{array}{|c|}
\hline
\phi_1(T)\\
\hline\\
(\alpha + 3)\tau+ 2\tau^2 + (4\alpha + 4)\tau^3\\
  2\tau + 2\tau^2 + (4\alpha + 2)\tau^3\\
 (4\alpha + 4)\tau + 2\tau^2 + 3\tau^3\\
 (4\alpha + 2)\tau + 2\tau^2 + (\alpha + 1)\tau^3\\
 3\tau + 2\tau^2 + (\alpha + 3)\tau^3\\
 (\alpha + 1)\tau + 2\tau^2 + 2\tau^3\\
 \hline
\end{array}
$
$
\begin{array}{|c|}
\hline
\phi_2(T)\\
\hline\\
(\alpha + 3)\tau + 2\tau^2 + 3\tau^3\\
 2\tau + 2\tau^2 + (\alpha + 1)\tau^3\\
 (4\alpha + 4)\tau + 2\tau^2 + (\alpha + 3)\tau^3\\
 (4\alpha + 2)\tau + 2\tau^2 + 2\tau^3\\
 3\tau + 2\tau^2 + (4\alpha + 4)\tau^3\\
 (\alpha + 1)\tau + 2\tau^2 + (4\alpha + 2)\tau^3\\
 \hline
\end{array}
$\\
$
\begin{array}{|c|}
\hline
\phi_3(T)\\
\hline\\
(\alpha + 3)\tau + (2\alpha + 1)\tau^2 +  (\alpha + 3)\tau^3\\
 2\tau + (2\alpha + 1)\tau^2 +  2\tau^3\\
 (4\alpha + 4)\tau + (2\alpha + 1)\tau^2 + (4\alpha + 4)\tau^3\\
 (4\alpha + 2)\tau + (2\alpha + 1)\tau^2 + (4\alpha + 2)\tau^3\\
 3\tau + (2\alpha + 1)\tau^2 +  3\tau^3\\
 (\alpha + 1)\tau + (2\alpha + 1)\tau^2 + (\alpha + 1)\tau^3\\
 \hline
\end{array}
$
$
\begin{array}{|c|}
\hline
\phi_4(T)\\
\hline\\
(\alpha + 3)\tau + 4\alpha\tau^2 + 2\tau^3\\
 2\tau + 4\alpha\tau^2 + (4\alpha + 4)\tau^3\\
 (4\alpha + 4)\tau + 4\alpha\tau^2 +  (4\alpha + 2)\tau^3\\
 (4\alpha + 2)\tau + 4\alpha\tau^2 + 3\tau^3\\
 3\tau + 4\alpha\tau^2 + (\alpha + 1)\tau^3\\
 (\alpha + 1)\tau + 4\alpha\tau^2 + (\alpha + 3)\tau^3\\
 \hline
\end{array}
$\\
$
\begin{array}{|c|}
\hline
\phi_5(T)\\
\hline\\
(\alpha + 3)\tau + 4\alpha\tau^2 + (\alpha + 1)\tau^3\\
 2\tau +  4\alpha\tau^2 + (\alpha + 3)\tau^3\\
 (4\alpha + 4)\tau +  4\alpha\tau^2 + 2\tau^3\\
 (4\alpha + 2)\tau + 4\alpha\tau^2 + (4\alpha + 4)\tau^3\\
 3\tau + 4\alpha\tau^2 + (4\alpha + 2)\tau^3\\
 (\alpha + 1)\tau + 4\alpha\tau^2 + 3\tau^3\\
 \hline
\end{array}
$
$
\begin{array}{|c|}
\hline
\phi_6(T)\\
\hline\\
(\alpha + 3)\tau + (3\alpha + 3)\tau^2 + (\alpha + 3)\tau^3\\
 2\tau + (3\alpha + 3)\tau^2 + 2\tau^3\\
 (4\alpha + 4)\tau + (3\alpha + 3)\tau^2 + (4\alpha + 4)\tau^3\\
 (4\alpha + 2)\tau + (3\alpha + 3)\tau^2 + (4\alpha + 2)\tau^3\\
 3\tau (3\alpha + 3)\tau^2 + 3\tau^3\\
 (\alpha + 1)\tau + (3\alpha + 3)\tau^2 + (\alpha + 1)\tau^3\\
 \hline
\end{array}
$\\
$
\begin{array}{|c|}
\hline
\phi_7(T)\\
\hline\\
(\alpha + 3)\tau + (\alpha + 4)\tau^2 + 2\tau^3\\
 2\tau + (\alpha + 4)\tau^2 + (4\alpha + 4)\tau^3\\
 (4\alpha + 4)\tau + (\alpha + 4)\tau^2 + (4\alpha + 2)\tau^3\\
 (4\alpha + 2)\tau + (\alpha + 4)\tau^2 + 3\tau^3\\
 3\tau + (\alpha + 4)\tau^2 + (\alpha + 1)\tau^3\\
 (\alpha + 1)\tau + (\alpha + 4)\tau^2 + (\alpha + 3)\tau^3\\
 \hline
\end{array}
$
$
\begin{array}{|c|}
\hline
\phi_8(T)\\
\hline\\
(\alpha + 3)\tau + (\alpha + 4)\tau^2 + (\alpha + 1)\tau^3\\
 2\tau + (\alpha + 4)\tau^2 + (\alpha + 3)\tau^3\\
 (4\alpha + 4)\tau + (\alpha + 4)\tau^2 + 2\tau^3\\
 (4\alpha + 2)\tau + (\alpha + 4)\tau^2 + (4\alpha + 4)\tau^3\\
 3\tau + (\alpha + 4)\tau^2 + (4\alpha + 2)\tau^3\\
 (\alpha + 1)\tau + (\alpha + 4)\tau^2 + 3\tau^3\\
\hline
\end{array}
$

\newpage
\bibliography{biblio}
\bibliographystyle{plain}
\end{document}